\documentclass[12pt]{amsart}
\usepackage{amssymb}

\newtheorem{theorem}{Theorem}
\newtheorem{proposition}{Proposition}
\newtheorem{conjecture}{Conjecture}
\newtheorem{definition}{Definition}
\newtheorem{lemma}{Lemma}

\begin{document}

\title{Complements of caustics of real function singularities}

\author{V.A.~Vassiliev}\address{Weizmann Institute of Science} 
\email{ vavassiliev@gmail.com}

\thanks{
This work was supported by the Absorption Center in Science of the Ministry of Immigration and Absorption of the State of Israel.}

\begin{abstract}
We study the topology of complements of caustics of function singularities of low codimensions, in particular \\
1) complete the enumeration of connected components of the complements of caustics of {\em simple} (in the sense of V.~Arnold) singularities, in particular find the numbers of these components for the last two classes, $E_7$ and $E_8$, remaining unknown after the works of R.~Thom, V.~Arnold  and V.~Sedykh; \\
2) realize all these components for simple singularities by explicitly constructed functions, and realize  their one-dimensional homology and cohomology groups by cycles and cocycles; \\
3) prove that (in contrast to the case of simple singularities) for some parabolic singularities the two-dimensional homology groups of the complements of their caustics are nontrivial.
\end{abstract}



\maketitle

\section{Introduction}

The caustics of ray systems appear, for example, as sunbeams, i.e., sets of concentration of rays, and hence as singular sets of luminosity. In terms of {\it generating families} of functions describing an optical state, it are the sets of parameter values corresponding to the functions with non-Morse critical points. In the case of non-degenerate families, they form hypersurfaces in the parameter spaces. The singularity types of caustics that characterize the concentration of rays coincide with the classes of the natural classification of critical points of functions, see \cite{A}, \cite{AVGZ}. The local geometry of a caustic close to its various points is determined by the types of these points with respect to this classification. In particular, a caustic can divide the neighborhood of such a point into many connected components, the number of which depends on the singularity class. The luminosity function is regular outside the caustic, but its behavior usually varies significantly in different components of its complement. 
 
 Therefore, the enumeration of these local components of complements of caustics is an early step in the study of luminosity functions at their singular points. 

A natural starting segment of the classification of function (and caustic) singularities consists of so-called simple singularities distinguished in \cite{A}. Namely, it consists of the classes $A_k$ ($k \geq 1$), $D_k$ ($k \geq 4$), $E_6,$ $ E_7$ and $E_8$. (In addition, some of these classes have several different real forms, see Table \ref{t1} below.)

For pictures of standard singularities of caustics of families depending on $\leq 3$ parameters (that is, of classes $A_2, A_3, A_4$ and $D_4$), see, for example, \cite{AVGZ}, \S 21.3 and \cite{AGLV2}, Chapter 2; apparently they were first sketched in \cite{thom} and appeared in a comprehensible form in \cite{Acongr}. V.D.~Sedykh \cite{Sed15}, \cite{selast} studied in detail the stratification of caustics at all other simple singularities, except for only two, namely at singularities of types $A_k$, $D_k$ and $E_6$. This study implies (by a kind of Alexander duality considerations) the homology groups of their local complements, in particular the numbers of local connected components of these complements, and the fact that their homology groups are trivial in all dimension greater than 1.

In the first part of the present work (sections 2 and 3) we find the numbers of such local components for the two remaining classes, $E_7$ and $E_8$, and explicitly implement these numbers for all simple singularities by functions representing all components. We also describe a simple invariant separating all these components. 

In \S \ref{onedim} we describe a one-dimensional cohomology class of complements of caustics that generates one-dimensional cohomology groups of all such components for singularities $A_k$, $D_k$ and $E_6$ (and conjecturally also for $E_7$ and $E_8$) for which this group is non-trivial. We also demonstrate the 1-cycles in function spaces that generate the homology groups of all these components (except maybe for $E_7$ and $E_8$).

In \S \ref{H2} we prove that the complements of the caustics of singularities of the classes $P_8^2$ and $J_{10}^3$ (which are some of the simplest non-simple singularities) have non-trivial two-dimensional homology groups.

\subsection{Main definition}
\label{maindef}

 For an introduction to the theory of caustics see, for example, \cite{AVGZ}, \cite{AGLV2}.

In short, let $f:( {\mathbb R}^n,0) \to ( {\mathbb R},0)$, $d f(0)=0$, be a smooth function with a critical point at the origin, and $F:( {\mathbb R}^n \times  {\mathbb R}^l,0) \to ( {\mathbb R},0)$ be its arbitrary smooth deformation, i.e. a family of functions $f_\lambda \equiv F(\cdot, \lambda):  {\mathbb R}^n \to  {\mathbb R},$ $\lambda \in  {\mathbb R}^l$, $f_0 \equiv f$. The {\it caustic} of this deformation is the set of values of the parameter $\lambda \in  {\mathbb R}^l$ such that the corresponding function $f_\lambda$ has at least one non-Morse critical point close to the origin in $ {\mathbb R}^n$. This is a subvariety in $ {\mathbb R}^l$ (of codimension 1 in all interesting cases) that usually splits a neighborhood of the origin of $ {\mathbb R}^l$ into several connected components. These components, their enumeration and topology are the main concern of this work.

For the parallel study of complements of some other bifurcation sets of deformations of simple functions, the so-called {\em discriminants} (that is,  sets of parameter values $\lambda$ such that the corresponding varieties $f_\lambda^{-1}(0) $ are singular), see, in particular, \cite{Lo}, \cite{sedykh}, \cite{Visr}. 

\subsection{Simple singularities and their deformations}

According to \cite{A}, any germ of a {\em simple} function singularity in $n\geq 2$ variables can be reduced by a choice of local coordinates to the form $\varphi(x, y)+Q(z_1, \dots, z_{n-2})$, where $\varphi$ is one of the polynomial normal forms listed in Table \ref{t1}, and $Q$ is a non-degenerate quadratic form. 
\begin{table}
\caption{Normal forms of real simple singularities in $ {\mathbb R}^2$ \qquad \qquad \ }
\label{t1}
\begin{tabular}{|l|l|c|c|}
\hline
Notation & Normal form & \# of components & $b_1$ \\ 
\hline 
$A_\mu$, \quad $\mu \ge 1$ & $\pm x^{\mu+1} \pm y^2 $ & $\left[\frac{\mu+2}{2}\right]$ & 0 \\ [3pt]
$D_{2k}^-$;\quad $k \ge 2$ & $x^2y \/- y^{2k-1} $ & $\frac{k^2+k}{2}$ & $\frac{k^2-k}{2}$ \\ [3pt]
$D_{2k}^+$;\quad $k \ge 2$ & $x^2y \/+ y^{2k-1} $ & $\frac{k^2+3k-2}{2}$ & $\frac{k^2-3k+2}{2}$ \\ [4pt]
$\pm D_{2k+1}$; \ $k \ge 2$ & $\pm (x^2y \/+ y^{2k})$ & $\frac{k^2+3k}{2}$ & $\frac{k^2-k}{2}$\\ [3pt]
$\pm E_6$ & $x^3 \pm y^4 $ & 7 & 1 \cr \hline
$E_7$ & $x^3 + x y^3 $ & 10 & $ \geq 5$ \cr
$E_8$ & $x^3 + y^5$ & 15 & $\geq 7$ \cr
\hline \end{tabular}
\end{table} 
The notations $A_\mu$, $D^-_{2k}$, etc., listed in its first column, are the names of {\it stable equivalence} classes, i.e. sets of all singular function germs reducible to each other by the change of local coordinates and the quadratic form $Q$ (including changes of the number $n-2$ of variables of $Q$, and hence also of the number $n$ of variables of entire function).
The local geometry of the caustics of all sufficiently large (so-called {\em versal}, see \cite{AVGZ}) deformations of stably equivalent singularities is the same,  so we can and will assume that all our functions depend on two variables only. 

We will consider {\em miniversal} deformations of these singularities, i.e. versal deformations depending on the minimal possible number of parameters, which is equal to the {\em Milnor number} $\mu(f)$ of the deformed function $f$ (these numbers are indicated by subscripts in the notation $A_\mu$, $D_{2k}^-$, etc. of singularity classes).

The canonical miniversal deformations of simple singularities are as follows:
\begin{eqnarray}
\ \ A_\mu & & f + \lambda_1 + \lambda_2x + \lambda_3x^2 + \dots + \lambda_{\mu} x^{\mu-1} \label{avd} \\
\ \ D_\mu & & f + \lambda_1 + \lambda_2x + \lambda_3y+
\lambda_4y^2 + \lambda_5 y^3 + \cdots + \lambda_{\mu}y^{\mu -2} \label{dvd} \\
\ \ E_6 & & f + \lambda_1 + \lambda_2x + \lambda_3y +
\lambda_4x y + \lambda_5y^2 + \lambda_6 x y^2 \label{e6vd} \\
\ \ E_7 & & f + \lambda_1 + \lambda_2x + \lambda_3y +
\lambda_4x y + \lambda_5y^2 + \lambda_6y^3 + \lambda_7y^4 \label{e7vd} \\ 
\ \ E_8 & & f + \lambda_1 + \lambda_2x + \lambda_3y +
\lambda_4x y + \lambda_5y^2 + \lambda_6x y^2 +
\lambda_7y^3 + \lambda_8x y^3 \label{e8vd} \ , 
\end{eqnarray}
where $f$ is the initial singularity given by the corresponding normal form of the Table \ref{t1}, and $\lambda_i$ are parameters of the deformation.

 Moreover, the caustic of any of these deformations is invariant under additions of constant functions, therefore we will also consider {\em shortened miniversal deformations} depending on $\mu(f)-1$ parameters and obtained from formulas (\ref{avd})--(\ref{e8vd}) by omitting the parameter $\lambda_1$. In what follows, we consider only the caustics in parameter spaces of canonical deformations (\ref{avd})--(\ref{e8vd}) of simple singularities or in the corresponding shortened deformations.

\subsection{Numbers of components of complements of caustics of simple singularities}

\begin{theorem}[V.D.~Sedykh, \cite{selast}]
\label{thm1}
The number of local components of the complement of the caustic of any versal deformation of a real function singularity of type $A_\mu$, $D_\mu$ or $E_6$ is given in the third column of Table \ref{t1}. All these components are either contractible or homotopy equivalent to circles; the numbers of components homotopy equivalent to circles for these singularities are given in the fourth column of Table \ref{t1}.
\end{theorem} 

In \S \ref{real} below we give some Morse perturbations representing all these components and do the same also for singularities $E_7$ and $E_8$.

\begin{theorem}
\label{them}
The number of local components of the complement of the caustic of a versal deformation of any real function singularity of type $E_7$ $($respectively, $E_8)$ is equal to 10 $($respectively, 15$)$. The number of components of the complement of the caustic of any versal deformation of a singularity of type $E_7$ $($respectively, $E_8)$, whose one-dimensional homology group contains the summand $ {\mathbb Z}$, is at least 5 $($respectively, 7$)$.
\end{theorem}

\begin{conjecture}
For singularities of classes $E_7$ and $E_8$

1. The $i$-dimensional homology groups of the complements of caustics of their versal deformations are trivial for any $i>1$,

2. The one-dimensional homology groups of these complements are isomorphic respectively to $ {\mathbb Z}^5$ and $ {\mathbb Z}^7$.
\end{conjecture}

Here is a partial result towards this conjecture.

\begin{proposition}
\label{part}
The complement of the caustic of the singularity $E_7$ (respectively, $E_8$)  contains exactly four (respectively, five) connected components consisting of morsifications, all whose critical points are real. All these components are homeomorphic to convex domains in ${\mathbb R}^7$ (respectively, ${\mathbb R}^8$).
\end{proposition}

\subsection{Basic invariant of components}
\label{basinv}

\begin{definition} \rm
The {\em passport} of a Morse function in two variables is the set 
$(m_-, m_0,$ $ m_+)$ of the numbers of its local minima, saddlepoints and maxima.
\end{definition}

It follows immediately from the definition that the passport is an invariant of the components of the complement of a caustic. The number $m_- + m_+ - m_0$ is the same for all morsifications of a given isolated singularity: it is equal to the local index of its gradient vector field. 

\begin{theorem} 
The passport is a complete invariant of the complement of the caustic of a simple singularity.
\end{theorem}

A proof of this theorem follows immediately from the Theorems \ref{thm1} and \ref{them} and the realizations of these components in \S \ref{real} below. Namely, 
for any simple singularity class, all the passports of its morsifications shown in \S \ref{real} are different, and the numbers of these morsifications are equal to the numbers of components given in these two theorems.

\subsection{Basic 1-cohomology class}
\label{bascoc}

For any component of the complement of a caustic, consisting of Morse functions with $m$ saddlepoints, consider the map from this component to the (unordered) configuration space $B( {\mathbb R}^2, m)$ taking any point $\lambda$ to the collection of saddlepoints of the function $f_\lambda$. 
If $m \geq 2$, then the group $H^1(B( {\mathbb R}^2, m),  {\mathbb Z})$ is equal to $ {\mathbb Z}$ and is generated by the {\em winding number}. Namely, identifying $ {\mathbb R}^2$ with the complex line, we can associate with any configuration $(z_1, \dots, z_m) \subset B( {\mathbb R}^2,m) \simeq B( {\mathbb C}^1, m)$ its  discriminant, i.e. the product of all $m(m-1)$ complex numbers $z_i - z_j$, $i \neq j$. The winding number is induced by this map from the canonical generator of the group $H^1( {\mathbb C}^*,  {\mathbb Z})$.

\begin{theorem}
 The winding numbers of configurations of saddlepoints generate the 1-dimensional cohomology groups of all components of complements of the caustics of all singularities of types $D_\mu$ and $E_6$.
The number of components of the complement of the caustic of any versal deformation of any singularity of type $E_7$ $($respectively, $E_8)$, on which the cohomology class of the winding number is non-trivial, is at least 5 $($respectively, 7$)$.
\end{theorem}

This theorem will be proved in \S \ref{onedim}. The last statement of Theorem \ref{them} follows immediately from it.

\subsection{Higher homology groups for non-simple singularities}

By Theorem \ref{thm1}, the complements of the caustics of singularities $A_\mu, D_\mu $ and $E_6$ do not have non-trivial homology groups in dimensions greater than 1. Conjecturally this is true also for the remaining two simple singularity classes, $E_7$ and $E_8$. On contrary, for some more complicated singularities this is not the case. 

Namely, the next in difficulty class of function singularities are the {\em parabolic} ones, see e.g. \cite{AVGZ}. There are eight stable equivalence classes of real parabolic singularities, $P_8^1, P_8^2$, $\pm X_9^0$, $X_9^1$, $X_9^2$, $J_{10}^1$ and $J_{10}^3$. We will show in \S \ref{H2} that for at least two of these classes, $P_8^2$ and $J_{10}^3$, the two-dimensional homology group of the complement of the caustic is non-trivial.

\section{Realization of components of complements of caustics}
\label{real}
\subsection{$A_\mu$} 
\label{aaa}
The case $A_\mu$ is trivial: the caustic of any singularity of this class in the parameter space $ {\mathbb R}^\mu$ of miniversal deformation (\ref{avd}) is diffeomorphic to the product of $ {\mathbb R}^1$ and the {\em discriminant} (see \S \ref{maindef}) of the miniversal deformation of singularity $A_{\mu-1}$; this diffeomorphism is provided by the differentiation with respect to $x$ and dilation of coordinates $\lambda_i$. In particular, the number of real critical points of a Morse perturbation $f_\lambda$ is a complete invariant of the component containing the parameter value $\lambda$; this invariant can take values $\mu$, $\mu -2$, $\mu-4...$ All these components can be represented by integrals of arbitrary polynomials in $x$ of degree $\mu$ with leading term $\frac{x^\mu}{\mu+1}$ and exactly $\mu, \mu-2, \dots$ simple roots.

\subsection{$D_{2k}^-$} 

\begin{figure}
\unitlength 0.80mm
\linethickness{0.4pt}
\begin{center}
\begin{picture}(119.00,25.00)
\put(58,9.3){\circle*{1.33}}
\put(58,18.5){\circle*{1.33}}
\put(7.00,14.00){\circle*{1.33}}
\put(27.00,14.00){\circle*{1.33}}
\put(47.00,14.00){\circle*{1.33}}
\put(58,2){\line(0,1){23}}
\put(53,13){$+$}
\put(60,13){$-$}
\bezier{120}(27.00,14.00)(17.00,3.00)(7.00,14.00)
\bezier{120}(27.00,14.00)(17.00,25.00)(7.00,14.00)
\bezier{52}(27.00,14.00)(30.00,19.00)(37.00,19.00)
\bezier{44}(37.00,19.00)(43.00,19.00)(46.50,14.50)
\bezier{44}(27.00,14.00)(30.00,9.00)(37.00,9.00)
\bezier{40}(37.00,9.00)(43.00,9.00)(46.50,13.50)
\bezier{60}(47.50,14.50)(57,23)(66.50,14.50)
\bezier{60}(47.50,13.50)(57,5)(66.50,13.50)
\put(6.50,14.50){\line(-1,1){5}}
\put(6.50,13.50){\line(-1,-1){5}}
\put(1,12.5){$+$}
\put(17.00,14.00){\makebox(0,0)[cc]{$+$}}
\put(37.00,14.00){\makebox(0,0)[cc]{$+$}}
\put(67.00,14.00){\circle*{1.33}}
\put(87.00,14.00){\circle*{1.33}}
\put(107.00,14.00){\circle*{1.33}}
\bezier{100}(67.50,14.50)(77.00,24.00)(86.50,14.50)
\bezier{100}(67.50,13.50)(77.00,4.00)(86.50,13.50)
\bezier{100}(87.50,14.50)(97.00,24.00)(106.50,14.50)
\bezier{100}(87.50,13.50)(97.00,4.00)(106.50,13.50)
\multiput(107.50,14.50)(0.12,0.12){50}{\line(0,1){0.12}}
\multiput(113.50,7.50)(-0.12,0.12){50}{\line(0,1){0.12}}
\put(111.00,14.00){\makebox(0,0)[cc]{$-$}}
\put(97.00,14.00){\makebox(0,0)[cc]{$-$}}
\put(77.00,14.00){\makebox(0,0)[cc]{$-$}}
\put(26,4){$-$}
\put(26,20){$-$}
\put(85,4){$+$}
\put(85,20){$+$}
\end{picture}
\end{center}
\caption{A morsification for $D_{2k}^-$ \ ($k=7$)}
\label{d89}
\end{figure}
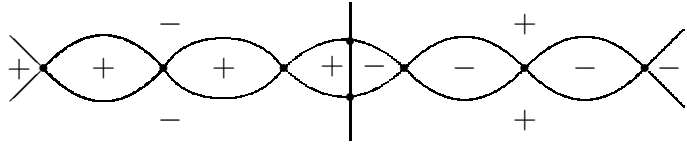
Any singularity of type $D_{2k}^-$ has a morsification, whose zero set is as shown in Fig.~\ref{d89}, where two non-vertical components of this set have $k-1$ intersection points and the vertical line separates these points into subsets of arbitrary cardinalities $r $ on the left and $k - 1 - r$ on the right, $r=0, 1, \dots, k-1$ (see e.g. \cite{Visr}). For any such $r$, the passport of corresponding morsifications is equal to $(r, k+1, k-1-r)$, so we get $k$ different passports.

Further, for any natural $ q= 2, \dots, k-1$ the preliminary perturbation $x^2 y - y^{2k-1} - \varepsilon y^{2q-1}$, $\varepsilon > 0$, has a singularity of type $D^-_{2q} $ at the origin and $2(k-q)$ non-real critical points. Applying to its real singularity all perturbations described in the previous paragraph, we obtain $q$ new morsifications with passports $(r, q+1, q-1-r)$, $r=0, 1, \dots, q-1$.
Finally, the initial function has the morsification $x^2 y - y^{2k-1} - \varepsilon y$, $\varepsilon >0$ with exactly two critical points and passport $(0, 2, 0)$. 
In total we have constructed $k + (k-1) + \dots +2 + 1= \frac{k(k+1)}{2}$ morsifications with different passports. By \cite{selast}, this is the number of all components of the complement of the caustic, so we have realized all these components (and proved that all of them have different passports).

\begin{figure}
\unitlength 0.70mm
\linethickness{0.4pt}
\begin{center}
{
\begin{picture}(81.00,27.00)
\put(21.00,15.00){\circle*{1.33}}
\put(41.00,15.00){\circle*{1.33}}
\bezier{100}(40.50,14.50)(31.00,5.00)(21.50,14.50)
\bezier{100}(21.50,15.50)(31.00,25.00)(40.50,15.50)
\bezier{100}(60.50,14.50)(51.00,5.00)(41.50,14.50)
\bezier{100}(41.50,15.50)(51.00,25.00)(60.50,15.50)
\bezier{44}(61.50,15.50)(64.00,19.00)(71.00,20.00)
\bezier{44}(61.50,14.50)(64.00,11.00)(71.00,10.00)
\bezier{60}(71.00,10.00)(81.00,10.00)(81.00,15.00)
\bezier{60}(71.00,20.00)(81.00,20.00)(81.00,15.00)
\bezier{44}(20.50,14.50)(18.00,11.00)(11.00,10.00)
\bezier{44}(20.50,15.50)(18.00,19.00)(11.00,20.00)
\put(61.00,15.00){\circle*{1.33}}
\put(10.00,20.00){\circle*{1.33}}
\put(10.00,10.00){\circle*{1.33}}
\bezier{52}(9.33,10.00)(1.00,10.00)(1.00,15.00)
\bezier{52}(1.00,15.00)(1.00,20.00)(9.33,20.00)
\put(10.00,27.00){\line(0,-1){6.33}}
\put(10.00,19.33){\line(0,-1){8.66}}
\put(10.00,9.33){\line(0,-1){6.33}}
\put(6.00,15.00){\makebox(0,0)[cc]{$+$}}
\put(14.00,15.00){\makebox(0,0)[cc]{$-$}}
\put(31.00,15.00){\makebox(0,0)[cc]{$-$}}
\put(51.00,15.00){\makebox(0,0)[cc]{$-$}}
\put(70,14){$-$}
\put(37,5){$+$}
\put(37,23){$+$}
\end{picture} \qquad \qquad \qquad
\begin{picture}(69.00,32.00)
\put(12.00,15.00){\circle*{1.33}}
\put(32.00,15.00){\circle*{1.33}}
\bezier{100}(31.50,15.50)(22.00,25.00)(12.50,15.50)
\bezier{100}(12.50,14.50)(22.00,5.00)(31.50,14.50)
\bezier{44}(32.50,14.50)(36.00,10.00)(42.00,10.00)
\bezier{40}(32.50,15.50)(36.00,20.00)(42.00,20.00)
\put(52.00,15.00){\circle*{1.33}}
\put(60,14){$-$}
\bezier{44}(52.50,14.50)(56,10)(62,10)
\bezier{40}(52.50,15.50)(56,20)(62,20)
\bezier{44}(51.50,14.50)(47,10)(42,10)
\bezier{40}(51.50,15.50)(47,20)(42,20)
\bezier{48}(62.00,20.00)(69.00,20.00)(69.00,15.00)
\bezier{48}(69.00,15.00)(69.00,10.00)(62.00,10.00)
\put(11.50,14.50){\line(-1,-1){6}}
\put(11.50,15.50){\line(-1,1){6}}
\put(41,5){\line(0,1){20}}
\put(41.00,10.00){\circle*{1.33}}
\put(41.00,20.00){\circle*{1.33}}
\put(7.00,15.00){\makebox(0,0)[cc]{$+$}}
\put(22.00,15.00){\makebox(0,0)[cc]{$+$}}
\put(37.00,15.00){\makebox(0,0)[cc]{$+$}}
\put(46.00,15.00){\makebox(0,0)[cc]{$-$}}
\put(30,5){$-$}
\put(30,20){$-$}
\put(50,5){$+$}
\put(50,20){$+$}
\end{picture}
}
\end{center}
\caption{Morsifications for $D_{2k}^+$ $(k=5)$ and $+D_{2k+1}$ $(k=4)$} 
\label{d8p9}
\end{figure}
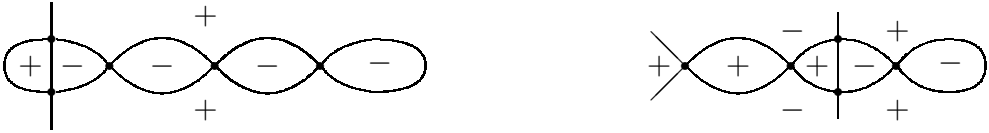

\subsection{$D_{2k}^+$} 
\label{d2kp}
Any singularity of type $D_{2k}^+$, $k \geq 2$, has morsifications with zero set as shown in Fig.~\ref{d8p9} (left), where the non-vertical component has $k-2$ self-intersection points, and the vertical component intersects it transversally and leaves arbitrarily many of these $k-2$ points to the left of it. If the intersection set of these components is non-empty, then the passport of the obtained morsification can be equal to either of numbers $(1, k, k-1), (2, k, k-2), \dots, (k-1, k, 1)$: in total $k-1$ possibilities. In addition, the vertical line can pass the other component of zero set from either side, which gives us two cases more with passports $(0, k-1, k-1)$ and $(k-1, k-1, 0)$. 

Further, for any $q = 2, \dots, k-1$ our singularity has the preliminary perturbation $x^2 y + y^{2k-1} + \varepsilon y^{2q-1}$ with a real singularity of class $D^+_{2q}$ and $2(k-q)$ non-real critical points. Applying to it the previous perturbations, we get another $q+1$ different passports. Finally, our singularity has the perturbation $x^2 y + y^{2k-1} + \varepsilon y$, $\varepsilon >0$, without real critical points (i.e. with passport $(0, 0, 0)$). In total we have realized $(k+1) + k + (k-1) + \dots +3 +1= \frac{k^2+ 3k -2}{2}$ different passports; by Theorem \ref{thm1} this is the number of components.

\subsection{$D_{2k+1}$} Any singularity of type $+D_{2k+1}$, $k \geq 2$, has morsifications with zero set as 
shown in Fig.~\ref{d8p9} (right), where the number of self-intersections of the non-vertical component is equal to $k-1$. Two components of this zero set can intersect each other in $k+1$ topologically different ways, providing $k+1$ different passports $(k, k+1, 0), (k-1, k+1, 1), \dots, (1, k+1, k-1)$, and finally $(0, k, k-1)$ for the case of empty intersection. Further, for any $q = 2, \dots, k-1$ our singularity has the preliminary perturbation $x^2 y + y^{2k} + \varepsilon y^{2q}$, $\varepsilon >0$, with a single singularity of type $+D_{2q+1}$ and $2(k-q)$ non-real critical points. Continuing all of these perturbations by above-described small morsifications of these their singularities, we get $(k+1) + k + \dots + 3 = \frac{(k+1)(k+2)}{2}-3$ different passports. In addition, there are Morse perturbations $x^2 y + y^{2k} + \varepsilon y^2 \pm \varepsilon y$ having passports $(1, 2, 0)$ or $(0, 1, 0)$ depending on the sign $\pm$. So, we have realized $\frac{k^2+3k}{2}$ different passports.

\subsection{$E_6$} 
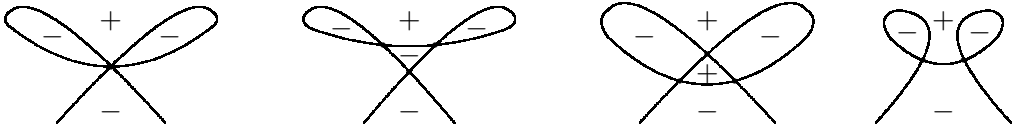
\begin{figure}
\unitlength 0.60mm
\linethickness{0.4pt}
\begin{center}
{
\begin{picture}(50,26)
\bezier{400}(13,1)(40,31)(47,26)
\bezier{100}(47,26)(50,24)(47,22)
\bezier{350}(47,22)(25,5)(3,22)
\bezier{100}(3,22)(0,24)(3,26)
\bezier{400}(3,26)(10,31)(37,1)
\put(25.00,13.5){\circle*{1.33}}
\put(38.00,20.00){\small \makebox(0,0)[cc]{$-$}}
\put(12.00,20.00){\small \makebox(0,0)[cc]{$-$}}
\put(22.5,2){\small $-$}
\put(22.5,22){\small $+$}
\end{picture}
\qquad 
\begin{picture}(50,26)
\bezier{400}(15,1)(40,31)(47,26)
\bezier{100}(47,26)(50,24)(47,22)
\bezier{320}(47,22)(25,14)(3,22)
\bezier{100}(3,22)(0,24)(3,26)
\bezier{400}(3,26)(10,31)(35,1)
\put(25.00,12.30){\circle*{1.33}}
\put(30.80,18.20){\circle*{1.33}}
\put(19.20,18.20){\circle*{1.33}}
\put(25.00,16.00){\small \makebox(0,0)[cc]{$-$}}
\put(40.00,22.00){\small \makebox(0,0)[cc]{$-$}}
\put(10.00,22.00){\small \makebox(0,0)[cc]{$-$}}
\put(22.5,2){\small $-$}
\put(22.5,22){\small $+$}
\end{picture}
\qquad 
\begin{picture}(50,26)
\put(25.00,16.20){\circle*{1.33}}
\put(31.10,10.20){\circle*{1.33}}
\put(18.90,10.20){\circle*{1.33}}
\bezier{400}(10,1)(40,34)(47,26)
\bezier{100}(47,26)(50,24)(47,20)
\bezier{400}(47,20)(25,-1)(3,20)
\bezier{100}(3,20)(0,24)(3,26)
\bezier{400}(3,26)(10,34)(40,1)
\put(25.00,12.00){\small \makebox(0,0)[cc]{$+$}}
\put(39.00,20.00){\small \makebox(0,0)[cc]{$-$}}
\put(11.00,20.00){\small \makebox(0,0)[cc]{$-$}}
\put(22.5,2){\small $-$}
\put(22.5,22){\small $+$}
\end{picture} 
\begin{picture}(50,26)
\put(29.5,15){\circle*{1.33}}
\put(20.5,15){\circle*{1.33}}
\bezier{400}(10,1)(31,24)(15,26)
\bezier{100}(13,20)(10,25)(15,26)
\bezier{400}(37,20)(25,8)(13,20)
\bezier{100}(35,26)(40,25)(37,20)
\bezier{400}(35,26)(19,24)(40,1)
\put(33.00,21.00){\small \makebox(0,0)[cc]{$-$}}
\put(17.00,21.00){\small \makebox(0,0)[cc]{$-$}}
\put(22.5,2){\small $-$}
\put(22.5,22){\small $+$}
\end{picture}
}
\end{center}
\caption{Perturbations for $+E_6$}
\label{e6proof}
\end{figure}

Singularity $+E_6$ has a perturbation with one singular point of type $D_4^-$ and two local minima, see the leftmost picture of Fig.~\ref{e6proof}. (A construction of the opposite perturbation of $-E_6$ singularity is described in p. 16 of \cite{AC}). Slightly deforming this singular point we obtain perturbations with passports $(3, 3, 0)$ and $(2, 3, 1)$ shown in the middle part of same figure. Also, the perturbation $x^3 + y^4 + \varepsilon y^3$ breaks the singularity $+E_6$ into two critical points of types $D_4^+$ and $A_2$. These singularities can be further perturbed independently, in particular the singularity $A_2$ can be eliminated (i.e. split into two non-real Morse points), and singularity $D_4^+$, according to \S \ref{d2kp}, can be moved to Morse functions with passports $(1, 2, 1)$, $(1, 1, 0)$, $(0, 1, 1)$ and $(0, 0, 0)$.

Further, the perturbation $x^3 + y^4 - 2 \varepsilon y^2 + \varepsilon^2$ splits singularity $+E_6$ into three points of type $A_2$, two of them with critical value 0 and the third with value $\varepsilon^2$. Let us perturb them independently in such a way that the point with value $\varepsilon^2$ vanishes, and each of two others splits into a saddlepoint with value 0 and a local minimum with a negative value. Zero level set of the obtained function looks as in the right-hand picture of Fig.~\ref{e6proof}, the passport of this function is $(2, 2, 0)$. So we have realized seven different passports, the maximal number allowed by Theorem \ref{thm1}. 

The case of singularity $-E_6$ can be reduced to that of $+E_6$ by the transformation $f_\lambda(x,y) \mapsto -f_\lambda(-x,y)$.

\unitlength 1mm
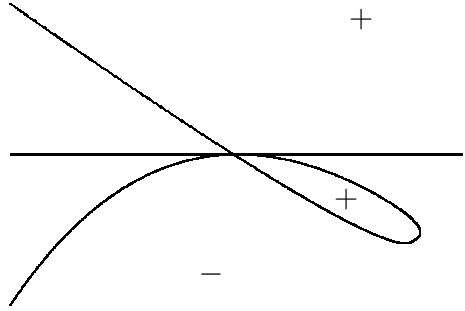
\begin{figure}
\begin{center}
\begin{picture}(60,40)
\put(0,20){\line(1,0){60}}
\bezier{500}(0,0)(20,30)(48,15)
\bezier{150}(48,15)(57,10)(53.5,8.5)
\bezier{500}(53.5,8.5)(50,5)(0,40)
\put(45,37){$+$}
\put(43,13){$+$}
\put(25,3){$-$}
\end{picture}
\end{center}
\caption{Perturbation $E_7 \mapsto D_6^-$}
\label{e7perper}
\end{figure}

\subsection{$E_7$} 
\label{e7p}
Zero level set of original singularity $x^3 + x y^3$ of type $E_7$ consists of two components: a semicubical parabola and its tangent line \
\unitlength 1.20mm
\linethickness{0.4pt}
 \begin{picture}(17,8)
\put(0,4){\line(1,0){15}}
\bezier{70}(0,0)(5,4)(13,4)
\bezier{70}(0,8)(5,4)(13,4)
\put(11,6){{\small $+$}}
\put(11,0){$-$}
\put(0,1.5){{\small $+$}}
\put(0,4.3){$-$}
\put(11.5,4){\circle*{0.8}}
\end{picture}. The perturbation $x^3 + x y^3 + \varepsilon x^2 y,$ $\varepsilon >0$, splits this singularity into a singular point of type $D_{6}^-$ and one Morse maximum, see Fig.~\ref{e7perper} (where the line $\{x=0\}$ is horizontal). Applying further to this singularity $D_6^-$ all perturbations described in \S \ref{d2kp}, we obtain morsifications with 
 passports $(0, 4, 3)$, $(1, 4, 2)$, $(2, 4, 1)$, $(0, 3, 2)$, $(1, 3, 1)$, and $(0, 2, 1)$. The space of deformation (\ref{e7vd}) admits the involution $f_\lambda(x, y) \leftrightarrow -f_\lambda(-x,y)$ which permutes functions with passports $(p, q, r)$ and $(r, q, p)$. Applying it to six morsifications obtained above, we get three new passports $(3,4,0)$, $(2,3,0)$, and $(1, 2, 0)$. Finally, the perturbation $x^3 + x y^3 + \varepsilon x y$, $\varepsilon >0$, has unique saddlepoint at the origin, i.e. realizes the tenth passport $(0,1,0)$. 

\subsection{$E_8$} 
\label{e8p}

 For any $\varepsilon \neq 0$, the perturbation $x^3 + y^5 + \varepsilon x y^3$ splits the singularity $E_8$ into a critical point of type $E_7$ and one point of local minimum or maximum depending on the sign of $\varepsilon$. Applying further to this singularity $E_7$ all perturbations described in \S \ref{e7p} we obtain functions with passports $(4, 4, 0)$, $(3, 4, 1)$, $(2, 4, 2)$, $(1, 4, 3)$, $(0, 4, 4)$,
$(3, 3, 0)$, $(2, 3, 1)$, $(1, 3, 2)$, $(0, 3, 3)$, $(2, 2, 0)$, $(1, 2, 1)$, $(0, 2, 2)$, $(1, 1, 0),$ and $(0, 1, 1)$. Finally, the perturbation $x^3 + y^5 + \varepsilon x$, $\varepsilon>0,$ has no real critical points, i.e. its passport is equal to $(0, 0, 0)$. In total we have realized 15 different passports.

\section{List of components is complete}
\label{nomore}

\begin{proposition}
\label{pronm}
Singularities $E_7$ and $E_8$ have no extra components of complements of caustics of versal deformations in addition to 10 and 15 components realized in \S\S \ref{e7p} and \ref{e8p}.
\end{proposition}

This proposition completes the proof of the first statement of Theorem \ref{them}. Our proof of this proposition follows the method of \S 6 of \cite{Visr}, where a similar statement on the complements of {\em discriminants} (i.e., sets of functions with zero critical value) was proved. In particular, it uses a computer algorithm enumerating all {\em virtual morsifications} of a singularity, i.e. collections of certain topological characteristics of its strict morsifications. 

Namely, such a set of characteristics of a non-discriminant real strict morsification $f_\lambda$ in $n$ real variables consists of
 \begin{enumerate}
\item the $\mu \times \mu$ matrix of intersection indices of canonically ordered and oriented vanishing cycles in the manifold $f^{-1}_\lambda(0)$, that are  defined by a standard system of paths connecting 0 with all critical values of $f_\lambda$, 
\item the string of $\mu$ intersection indices of these vanishing cycles with the cycle $f^{-1}_\lambda(0) \cap  {\mathbb R}^n$,
\item Morse indices of all real critical points of $f_\lambda$ ordered by increase of their critical values, and 
\item the number of these critical points having negative critical values.
\end{enumerate}

An open subset of the parameter space $ {\mathbb R}^\mu$ of any miniversal deformation of a real function singularity consists of non-discriminant strict morsifications. This subset splits into
 many open domains (non necessarily connected) consisting of such morsifications characterized by different virtual morsifications. These domains are separated by walls consisting of points of different bifurcation sets: caustic, discriminant, and also real and complex Maxwell sets (that is, closures of sets of parameter values corresponding to functions with coinciding critical values at different critical points). In the case of simple singularities, the data of a virtual morsification imply a complete information on the set of the walls bounding the corresponding domain of $ {\mathbb R}^\mu$, and on the neighboring open domains, to which one can pass from it by crossing a single wall; these passages are described by {\em virtual surgeries}, i.e. certain flips of virtual morsifications. Applying all possible sequences of these flips to a single virtual morsification, one can list all virtual morsifications related with real morsifications of our singularity.

The set of all virtual morsifications splits into {\em virtual components of complements of caustics}, i.e. sets of virtual morsifications which can be obtained from one another by chains of formal moves modeling the surgeries of real morsifications which do not cross the caustic. 
In particular, if two morsifications of a singularity belong to the same component of the complement of the caustic, then related virtual morsifications belong to the the same virtual component.

The passport of a real morsification can be read from the corresponding virtual morsification. Virtual morsifications from the same virtual component obviously have equal passports.

Investigating a singularity, we proceed as follows.
First we find a virtual morsification of our singularity describing a real morsification and substitute it as initial data to the non-restricted program, which enumerates all virtual morsifications and gives us the numbers of virtual morsifications with each admissible passport. This program (with initial data of $E_7$ singularity) can be found at
{\footnotesize
https://drive.google.com/drive/folders/1lsSxpValbfeXk51MD4D9rXpA3dk9zeAr?usp=sharing, 
}
file mE7gen. Further, we run the restricted program (see file mE7part) which counts all virtual morsifications within some virtual component (and fix the number of them). Then choose an arbitrary morsification from the first list but not from the second, substitute it to the restricted 
program, and get some other virtual component, etc. until the sum of numbers of virtual morsifications in the explored virtual components becomes equal to the total number provided by the non-restricted program.  

More detailed description of our algorithms can be read in \cite{APLT}, Chapter V, and \cite{Visr} with one modification: in these works the virtual components of complements of discriminants were considered, therefore in the reduced program we switch out another set of flips: the ones modeling the surgeries crossing the discriminant. 

The results of the work of this algorithm for singularities $E_7$ and $E_8$ are summarized in the next two propositions.

\begin{proposition}
\label{pro1}
The singularity $E_7$ has exactly 8648 different virtual morsifications split into 10 virtual components of the complement of the caustic. If two virtual morsifications of this list have equal passports, then they belong to the same virtual component. The numbers of virtual morsifications in these components are
 described by Table \ref{t2}, in which passports are characterized by only two numbers: the total number $M=m_- + m_0 + m_+$ of real critical points and the number $m_-$ of minima $($two other numbers can be deduced from them by equalities $m_0= (M+1)/2$ and $m_+ = m_0 -1 - m_-)$. \hfill $\Box$
\end{proposition}

\begin{table}
\caption{Statistics of virtual morsifications for $E_7$}
\label{t2}
\begin{tabular}{|c|cccccc|}
\hline
Passport data $(M,m_-)$ & (1,0)& (3,0) & (5,0) & (5,1) & (7,0) & (7,1) \\
 & & (3,1) & (5,2) & & (7,3) & (7,2) \\
\hline 
Number of virtual morsifications & 260 & 312 & 636 & 348 & 2384 & 688 \\
\hline
\end{tabular}
\end{table}

\begin{proposition}
\label{pro2}
The singularity $E_8$ has exactly 51468 different virtual morsifications split into 15 virtual components of complements of caustics. If two virtual morsifications of this list have equal passports, then they belong to the same virtual component. The numbers of virtual morsifications in these virtual components are
 described by Table \ref{t3}. \hfill $\Box$
\end{proposition}
\begin{table}
\caption{Statistics of virtual morsifications for $E_8$}
\label{t3}
\begin{tabular}{|c|ccccccccc|}
\hline
$(M,m_-)$ & (0,0)& (2,0) & (4,0) & (4,1) & (6,0) & (6,1) & (8,0) & (8,1) & (8,2) \\
& & (2,1) & (4,2) & & (6,3) & (6,2) & (8,4) & (8,3) & \\
\hline 
Number & 1200 & 819 & 1195 & 790 & 3227 & 1246 & 13599 & 3690 & 1926 \\
\hline
\end{tabular}
\end{table}

\noindent
{\bf Remark.}
In the case of arbitrary simple singularities any two virtual morsifications with equal passports belong to the same virtual component. This property fails already for the simplest non-simple singularities of classes $X_9$.
\medskip

It remains to prove that each virtual component of the complement of the caustic of singularity $E_7$ or $E_8$ mentioned in Propositions \ref{pro1} and \ref{pro2} is related with unique component of the complement of the real caustic in the corresponding parameter space $ {\mathbb R}^7$ or $ {\mathbb R}^8$.

The main idea behind our algorithm is that of the {\em Looijenga map} (see \cite{Lo0}) that acts from the parameter space $ {\mathbb C}^\mu$ of the complexified miniversal deformation of a singularity $f$ to the space (also equal to $ {\mathbb C}^\mu$) of polynomials of degree $\mu$ with top coefficient 1. This map associates any parameter value $\lambda \in  {\mathbb C}^\mu$ with the set of coefficients of the polynomial, whose roots are the critical values of the corresponding perturbation $f_\lambda$. 
In the case of simple complex singularities and their canonical deformations (\ref{avd})--(\ref{e8vd}) this map is proper, and its restriction to the set of strictly Morse functions is a covering over the set of polynomials with only simple roots. Therefore any movement of $\mu$ critical values of a strictly Morse function $f_\lambda$ can be lifted to a movement of functions, starting with $f_\lambda$ and having the prescribed critical values.
The restriction of the Looijenga map to the parameter space $ {\mathbb R}^\mu$ of the real versal deformation obviously maps it into the space of real polynomials.

\begin{lemma}
\label{pro3}
If two different components of the complement of the caustic in the parameter space $ {\mathbb R}^\mu$ of the canonical real versal deformation $($\ref{avd}$)$--$($\ref{e8vd}$)$ of a simple singularity contain strict morsifications characterized by the same virtual morsification, then there exists a linear automorphism of $ {\mathbb R}^\mu$ commuting with the Looijenga map and taking these components to one another. 
\end{lemma}

\noindent
{\it Proof.} The sets of all complex critical values of such two morsifications can be continuously deformed to each other avoiding collisions and respecting the complex conjugation. Therefore, using the Looijenga covering, we can assume that these two morsifications $f_\lambda, f_{\lambda'}$ have the same sets of critical values and the same system of paths defining the corresponding vanishing cycles. The situation when two different (generally, complex) morsifications of a simple singularity have the same critical values and the same intersection matrix of corresponding vanishing cycles defined by the same system of paths was studied in \cite{Liv}. It follows easily from Picard--Lefschetz formula that in this situation there exists an automorphism of the Looijenga covering on the space of strictly Morse complex functions of the corresponding canonical deformation (\ref{avd})--(\ref{e8vd}), that sends $\lambda$ to $\lambda'$. The group of all automorphisms of this covering for each simple singularity was calculated in \cite{Liv}: it turns out that all these automorphisms can be extended to linear automorphisms of entire parameter space $ {\mathbb C}^\mu$, and moreover they are induced by certain linear automorphisms of the argument space $ {\mathbb C}^n$ of our functions that define symmetries of the initial singularity $f$. 

If functions $f_\lambda$ and $f_{\lambda'}$ are real, then the automorphisms of $ {\mathbb C}^\mu$ and $ {\mathbb C}^n$ defined in this way preserve respectively the subspaces $ {\mathbb R}^\mu \subset  {\mathbb C}^\mu$ of real morsifications and $ {\mathbb R}^n \subset  {\mathbb C}^n$ of real arguments of our functions. \hfill $\Box$ \medskip

To prove Proposition \ref{pronm}, it remains to show that singularities $E_7$ and $E_8$ do not admit such non-trivial automorphisms of parameter spaces of corresponding real versal deformations (\ref{e7vd}), (\ref{e8vd}); we can assume that $n=2$.

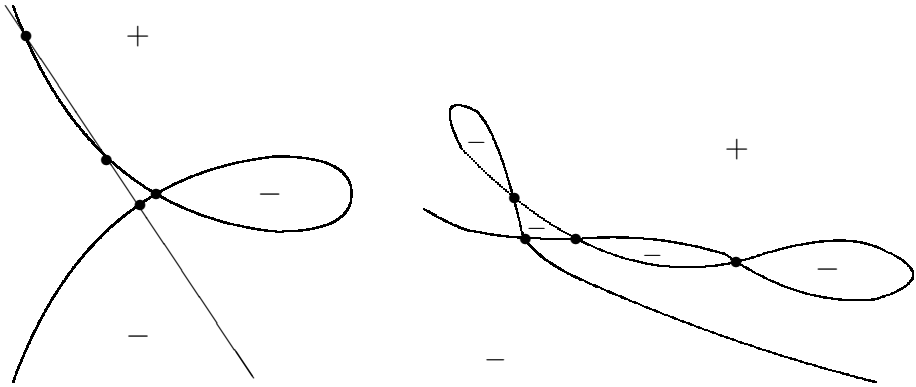
\begin{figure}
\unitlength 0,5mm
\begin{center}
\begin{picture}(100,100)
\bezier{800}(10,0)(30,55)(80,60)
\bezier{300}(80,60)(100,60)(100,50)
\bezier{300}(100,50)(100,40)(80,40)
\bezier{800}(10,100)(30,45)(80,40)
\put(8,100){\line(2,-3){66}}
\put(13.5,92){\circle*{3}}
\put(34.8,59){\circle*{3}}
\put(43.8,47){\circle*{3}}
\put(48,50){\circle*{3}}
\put(75,48){$-$}
\put(40,10){$-$}
\put(40,90){$+$}
\end{picture} \qquad 
\unitlength 1mm
\begin{picture}(65,36)
\bezier{70}(0,23)(5,20)(7,20)
\bezier{70}(7,20)(11,19)(20.2,19)
\bezier{100}(20.2,19)(30,20)(40,17)
\bezier{100}(40,17)(50,10)(60,11)
\bezier{70}(60,11)(70,14)(60,18)
\bezier{100}(60,18)(55,20)(45,17)
\bezier{150}(45,17)(25,10)(5,31)
\bezier{100}(5,31)(1,39)(7,36)
\bezier{100}(7,36)(10,33)(13,20)
\bezier{100}(13,20)(14,17)(20,14)
\bezier{150}(20,14)(40,5)(60,0)
\put(12,24.5){\circle*{1.5}}
\put(13.5,19){\circle*{1.5}}
\put(20.2,19){\circle*{1.5}}
\put(41.5,16){\circle*{1.5}}
\put(8,2){\small $-$}
\put(5.6,31){\footnotesize $-$}
\put(13.6,19.5){\footnotesize $-$}
\put(29,16){\footnotesize $-$}
\put(52,14){$-$}
\put(40,30){$+$}
\end{picture}
\end{center}
\caption{Nice perturbations for $E_7$ and $E_8$}
\label{e8proof}
\end{figure}

Singularities $E_7$ and $E_8$ have perturbations whose zero sets look topologically as shown respectively in the left and right parts of Fig.~\ref{e8proof}, see \cite{AC}, pages 16 and 17. (Moreover, the left-hand picture can be realized by easy formula $(x^2 + y^3 - \varepsilon^2 y^2)(x + 2 \varepsilon y + \frac{\varepsilon^3}{27}). ) $
Adding a small positive constant to either of these functions we obtain a morsification $f_\lambda$ of singularity $E_7$ (respectively, $E_8$) having only real critical points, namely 3 (respectively, 4) local Morse minima with negative critical values and 4 saddlepoints with positive values. It is easy to see that the entire component of the complement of the {\em discriminant} of our deformation containing this morsification consists of functions with the same property. Indeed, to change the number of critical points or their Morse indices, two critical points with different Morse indices should meet, but their critical values 
have different signs and cannot meet avoiding the value 0 which is prohibited for non-discriminant functions. 

By Theorem (1.7) of \cite{Lo} this component is contractible, in particular for any its point $\lambda$ any ordering of critical points of $f_\lambda$ can be uniquely extended by continuity to orderings of critical points of functions $f_{\lambda'}$ for all $\lambda'$ from this component. The ordered analog of the Looijenga map, sending any point $\lambda$ to the (ordered) collection of critical values at these ordered critical points, is a diffeomorphism of our component to an open octant in $ {\mathbb R}^\mu$. 

It is easy to check that the complex conjugation acts trivially on the homology group of the {\em Milnor fiber} (i.e. zero level set) of the stabilization $f_\lambda(x, y) + z^2:  {\mathbb C}^3 \to  {\mathbb C}$ of any function $f_{\lambda'}$ from this component. By Theorem (1.6) of \cite{Lo} the complement of the discriminant of singularity $E_7$ (respectively, $E_8$) admits only one connected component, whose points $\lambda'$ satisfy this property. In particular, any automorphism of $ {\mathbb R}^\mu$ (and hence $ {\mathbb C}^\mu$) commuting with Looijenga map as in the proof of Lemma \ref{pro3} sends this component to itself. Let us prove that this automorphism can be only the identical map.
For any function $f_\lambda$ from this component, the bilinear intersection form of vanishing cycles in 2-dimensional homology group of the Milnor fiber of function $f_\lambda(x, y) + z^2$ can be computed by the Gusein-Zade--A'Campo algorithm \cite{AC}, \cite{GZ}. This computation in the basis defined by a standard system of paths (see \cite{GZ}) shows that this form is defined by the standard (tree) Coxeter-Dynkin diagram of $E_7$ or $E_8$. 

 If $\lambda$ is a point of this component corresponding to a strictly Morse function $f_\lambda$, and $\lambda'$ is its image under an automorphism of $ {\mathbb R}^\mu$ commuting with Looijenga map, then we get two one-to-one correspondences between critical points of $f_\lambda$ and these of $f_{\lambda'}$: one of them preserves critical values, and the other is defined by the continuation over the paths in $ {\mathbb R}^\mu$ connecting $\lambda$ and $\lambda'$ and not intersecting the discriminant. Composing these two correspondences, we obtain a permutation of critical points of $f_\lambda$. By the construction, this permutation preserves the intersection indices of corresponding vanishing cycles, and hence defines an automorphism of the canonical Coxeter-Dynkin diagram of singularity $E_7$ or $E_8$. But these diagrams have no non-trivial automorphisms, therefore any desired automorphism or $ {\mathbb R}^\mu$ acts trivially on our component of the complement of the discriminant, and hence also on entire $ {\mathbb R}^\mu$. \hfill $\Box$

\subsection*{Proof of Proposition \ref{part}}
Let $f_\lambda$ be a strict morsification of a function singularity in two variables, all whose $\mu(f)$ critical points are real and critical values are not equal to 0. Consider a system of non-intersecting paths connecting the noncritical value 0 of this function with these critical values and the basis of vanishing cycles in the 2-dimensional homology group of the Milnor fiber of the function $f_\lambda(x,y)+z^2$ defined by these paths. 
We can assume that an orientation of the argument space ${\mathbb R}^3$ is chosen, then also some canonical orientations of these vanishing cycles are uniquely determined, see, for example, \cite{APLT}, \S V.1.6(iii). 

Consider the  Coxeter-Dynkin graph describing the intersection indices of these vanishing cycles, and orient  each of its edges from the vertex corresponding to the critical point with the higher critical value to that with the lower one. 
The isomorphism class of this oriented graph is the same for all functions $f_\lambda$ from the same connected component of the complement of the caustic. Indeed, a generic path inside the complement of the caustic can only reorder in ${\mathbb R}^1$ the neighboring critical values of $f_\lambda$ such that corresponding vertices are not connected by any edges of the graph, or move the critical values through 0. 

\begin{lemma}
For any connected component of the complement of the caustic of singularity $E_7$ or $E_8$ consisting of morsifications with only real critical points, the Coxeter-Dynkin diagram just defined has no non-trivial automorphisms. 
\end{lemma} 

\noindent
This follows immediately from the above listing of these components. \hfill $\Box$ \medskip

For this reason, the critical points of any such morsification can be ordered in a way that depends continuously of the morsification, and the ordered Looijenga map from this component to ${\mathbb R}^\mu$ ($\mu = 7$ or $8$) is well defined. By the properness of the Looijenga map, it is a proper submersion (and, therefore, a covering) from this connected component to the open domain in ${\mathbb R}^\mu$ distinguished by only the following restrictions (and their consequences): if our oriented Coxeter-Dynkin graph contains a segment directed from the $i$th vertex to the $j$th vertex, then the coordinate $w_i$ of the allowed point in ${\mathbb R}^\mu$ must be greater than $w_j$. This is a system of linear inequalities that defines a convex subset of ${\mathbb R}^\mu$. \hfill $\Box$

\section{One-dimensional homology of complements of caustics of simple singularities}
\label{onedim}

In this section we for any simple singularity class realize \ $b_1$ \ one-dimensional cycles in the complement of its caustic, where \ $b_1$ \ is the corresponding number written in the last column of Table \ref{t1}. All these cycles belong to different components of this complement, and the cocycle of \S \ref{bascoc} takes value 1 on each of them.

\subsection{Basic example} 
\label{basic}

Any function singularity of class $D_4^-$ in two variables  has the form $3x^2 y - y^3$ in appropriate local coordinates. Choose arbitrarily a number $\varepsilon \neq 0$ and consider the following one-parametric family of perturbations of this singularity:
\begin{equation}
f_\tau (x, y) \equiv 3x^2 y - y^3 - 3\varepsilon^2(x \sin \tau + y \cos \tau ), \quad \tau \in [0, 2\pi].
\label{param}
\end{equation}
It is easy to see that for any $\tau$ this function $f_\tau$ has exactly two real critical points, namely the saddlepoints whose coordinates $(x, y)$ are equal to \begin{equation}
\label{resad}
\pm \varepsilon \left(\cos \frac{\tau}{2}, \, \sin\frac{\tau}{2}\right),
\end{equation}
and two imaginary points with coordinates
\begin{equation}
\label{imsad}
\pm i \varepsilon \left(-\sin \frac{\tau}{2}, \, \cos \frac{\tau}{2}\right).
\end{equation}
 In particular, the family of functions (\ref{param}) permutes two saddlepoints (\ref{resad}), and the winding number of the corresponding loop in the configuration space $B( {\mathbb R}^2, 2)$ is equal to 1. 

\unitlength 1.00mm
\linethickness{0.4pt}
\begin{figure}
\begin{center}
\begin{picture}(40.00,58.00)
\put(22.00,27.00){\circle{1.00}}
\bezier{84}(22.30,27.40)(27.00,35.00)(23.00,46.00)
\bezier{124}(22.30,27.40)(29.00,38.00)(31.00,55.00)
\bezier{96}(22.30,27.40)(29.00,38.00)(38.00,45.00)
\bezier{56}(23.00,46.00)(28.00,47.00)(31.00,55.00)
\bezier{56}(38.00,45.00)(33.00,47.00)(31.00,55.00)
\bezier{84}(21.70,26.60)(19.00,20.00)(7.00,15.00)
\bezier{84}(7.00,15.00)(11.00,8.00)(24.00,7.00)
\bezier{80}(24.00,7.00)(13.00,7.00)(5.00,4.00)
\bezier{48}(5.00,4.00)(7.00,7.00)(7.00,15.00)
\bezier{36}(5.00,4.00)(9.00,6.00)(12.50,9.60)
\bezier{96}(24.00,7.00)(17.00,14.00)(21.70,26.60)
\bezier{12}(25.00,46.00)(27.00,46.00)(28.00,45.90)
\bezier{16}(31.00,45.85)(32.00,45.84)(35.50,45.50)
\bezier{30}(13.7,11.3)(17,15)(19,20)
\put(1.00,1.00){\makebox(0,0)[cc]{$A_3$}}
\put(5.00,19.00){\makebox(0,0)[cc]{$A_3$}}
\put(28.00,6.00){\makebox(0,0)[cc]{$A_3$}}
\put(26.00,26.00){\makebox(0,0)[cc]{$D_4^-$}}
\put(20.00,46.00){\makebox(0,0)[cc]{$A_3$}}
\put(31.00,58.00){\makebox(0,0)[cc]{$A_3$}}
\put(40.00,47.00){\makebox(0,0)[cc]{$A_3$}}
\put(14.00,3.00){\makebox(0,0)[cc]{$A_2$}}
\end{picture}
\caption{Pyramid: caustic for singularity $D_4^-$}
\label{pyramid}
\end{center}
\end{figure}
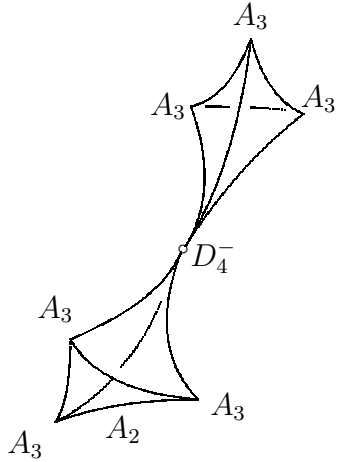

In a neighborhood of the origin of the parameter space $ {\mathbb R}^3$ of the shortened (i.e. with omitted parameter $\lambda_1$) versal deformation (\ref{dvd}) of singularity $D_4^-$, the caustic looks as shown in Fig.~\ref{pyramid}, see e.g. \cite{AVGZ}, \S 21. The loop (\ref{param}) hugs it at the waist. The interiors of two pyramids in this picture consist of morsifications with passports $(1, 3, 0)$ and $(0, 3, 1)$, and the exterior part consists of morsifications with passport $(0, 2, 0)$.

Now let $f$ be a function singularity, and $\lambda$ be a point of the parameter space $ {\mathbb R}^{\mu(f)}$ of its versal deformation such that the corresponding function $f_\lambda$ has one critical point of type $D_4^-$, and all other its real critical points are Morse. In a neighborhood of such a point $\lambda$, the pair $( {\mathbb R}^{\mu(f)}, \mbox{the caustic})$ is ambient diffeomorphic to the direct product of $ {\mathbb R}^{\mu(f)-3}$ and the pair consisting of the space $ {\mathbb R}^3$ and the standard caustic of $D_4^-$ shown in Fig.~\ref{pyramid}. The product of the origin of $ {\mathbb R}^{\mu -3}$ and our loop (\ref{param}) with sufficiently small $\varepsilon$ in $ {\mathbb R}^3$ embracing the caustic  defines then a 1-cycle in the complement of the caustic of $f$. If $f_\lambda$ has $p$ Morse minima, $q$ saddlepoints and $r$ maxima, then this loop will lie in the component of this complement with passport $(p, q+2, r)$. Our basic 1-cohomology class defined by winding numbers (see \S \ref{bascoc}) takes value 1 on this loop, since all other critical points of functions from this loop remain in small non-intersecting domains.

\subsection{1-cycles for simple singularities}

By Theorem 1.7 of \cite{Lo} all components of the complements of {\em discriminants} of all singularities $A_\mu$ are contractible, therefore the same is true for complements of caustics of $A_{\mu+1}$, cf. \S \ref{aaa}. \medskip

Any singularity of class $D_{2k}^-$ has perturbations, whose zero level sets look similar to Fig.~\ref{d89}, 
with two non-vertical components having an arbitrary number $r=1, 2, \dots, k-1$ of intersection points, but (unlike Fig.~\ref{d89}) the vertical line intersecting them at either of these $r$ intersection points. All these perturbations satisfy the conditions described in the previous subsection, i.e. have one point of type $D_4^-$ and only Morse remaining critical points.
Therefore our construction gives us loops in the components with all passports $(p, q, r)$ such that $q \in [2, k]$, $p \in [0, q-2]$ and $p+r=q-2$. There are exactly $\frac{k(k-1)}{2}$ such passports. \medskip

In the case of singularity $D_{2k}^+$, we proceed \label{1d2p} in a similar way with Fig.~\ref{d8p9} (left) and obtain loops in all components with passports $(p , q, r),$ where 
$p \geq 1, r \geq 1$, $q=p+r$ and $q \leq k-1$. There are exactly $\frac{(k-1)(k-2)}{2}$ such passports.

\medskip
In the case $+D_{2k+1}$, we get analogously the loops in all components with passports $(p, q, r)$ where $p+r = q-1$, $q \leq k$, $p \geq 1$. There are exactly $\frac{k(k-1)}{2}$ such passports.

\medskip
In the case $+E_6$ the needed perturbation is shown in the leftmost picture of Fig.~\ref{e6proof}, it gives us a 1-cycle in the component with passport $(2, 2, 0)$. 
\medskip

\unitlength 0.30mm
\begin{figure}
\begin{picture}(100,105)
\bezier{800}(10,0)(30,55)(80,60)
\bezier{300}(80,60)(94,60)(95,50)
\bezier{300}(95,50)(94,40)(80,40)
\bezier{800}(10,100)(30,45)(80,40)
\put(0,2){\line(1,1){98}}
\put(22,24){\circle*{3}}
\put(48,50){\circle*{4}}
\put(40,8){$-$}
\put(40,90){$+$}
\put(71,47){$+$}
\end{picture} \qquad
\begin{picture}(100,100)
\bezier{800}(10,0)(30,55)(80,60)
\bezier{300}(80,60)(94,60)(95,50)
\bezier{300}(95,50)(94,40)(80,40)
\bezier{800}(10,100)(30,45)(80,40)
\put(0,50){\line(1,0){100}}
\put(95,50){\circle*{3}}
\put(48,50){\circle*{4}}
\put(40,8){$-$}
\put(40,90){$+$}
\put(70,51){\small $-$}
\put(70,42){\small $+$}
\end{picture} \qquad 
\begin{picture}(100,100)
\bezier{800}(10,0)(30,55)(80,60)
\bezier{300}(80,60)(94,60)(95,50)
\bezier{300}(95,50)(94,40)(80,40)
\bezier{800}(10,100)(30,45)(80,40)
\put(0,98){\line(1,-1){98}}
\put(22,76){\circle*{3}}
\put(48,50){\circle*{4}}
\put(40,8){$-$}
\put(40,90){$+$}
\put(71,47){$-$}
\end{picture}
\caption{Degenerate perturbations for $E_7$}
\label{e7nfd}
\end{figure}
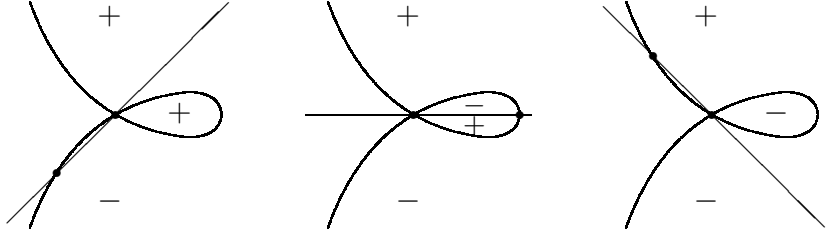

In the case $E_7$ we have three perturbations 
$$(x^2 + y^3 - \varepsilon^2 y^2)(x-2 \varepsilon y), \ (x^2 + y^3 - \varepsilon^2 y^2)x, \ (x^2 + y^3 - \varepsilon^2 y^2)(x + 2\varepsilon y) $$
shown in Fig.~\ref{e7nfd}. Our construction applied to their $D_4^-$-singularities realizes non-trivial 1-cycles in three components with passports $(2, 3, 0)$, $(1, 3, 1)$ and $(0, 3, 2)$.

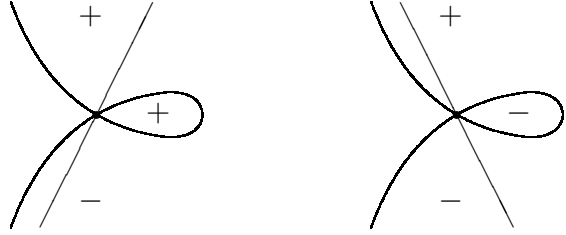
\begin{figure}
\begin{picture}(100,105)
\bezier{800}(10,0)(30,55)(80,60)
\bezier{300}(80,60)(94,60)(95,50)
\bezier{300}(95,50)(94,40)(80,40)
\bezier{800}(10,100)(30,45)(80,40)
\put(23,0){\line(1,2){50}}
\put(48,50){\circle*{4}}
\put(40,8){$-$}
\put(40,90){$+$}
\put(70,47){$+$}
\end{picture}
\qquad \qquad
\begin{picture}(100,105)
\bezier{800}(10,0)(30,55)(80,60)
\bezier{300}(80,60)(94,60)(95,50)
\bezier{300}(95,50)(94,40)(80,40)
\bezier{800}(10,100)(30,45)(80,40)
\put(23,100){\line(1,-2){50}}
\put(48,50){\circle*{4}}
\put(40,8){$-$}
\put(40,90){$+$}
\put(70,47){$-$}
\end{picture}
\caption{Two perturbations more for $E_7$}
\label{e7zh}
\end{figure}

Consider the stratum of the caustic of $E_7$ singularity containing the function $f_{\lambda_0} \equiv (x^2 + y^3 - \varepsilon^2 y^2)(x-2 \varepsilon y) $, whose zero set is shown in the leftmost picture of Fig.~\ref{e7nfd}. It consists of functions with one singularity of type $D_4^-$ and three real Morse points: two local maxima and one saddlepoint. By the properness of the Looijenga map of simple singularities (see \cite{Lo0}, \cite{APLT}), there exists a path $\{\lambda_t\},$ $t \in [0, 1),$ in this stratum, starting at this point $\lambda_0$ and such that all its points $f_{\lambda(t)}$ are functions having a point of type $D_4^-$ with zero critical value, the critical value at the right-hand maximum
points remaining constant, and the critical values at two other Morse points tending to one another and finally meeting at some intermediate value. The limit point of this path has a critical point of type $A_2$. By a small perturbation we can split this point into two non-real critical points, preserving the singularity of type $D_4^-$.
Zero level set of the obtained function looks topologically as shown in Fig.~\ref{e7zh} (left).

Applying the construction of \S \ref{basic} to its point of type $D_4^-$, we obtain a loop in the component with passport $(0, 2, 1)$. The symmetric construction applied to the rightmost picture of Fig.~\ref{e7nfd} gives us a cycle in the component with passport $(1, 2, 0)$, see Fig.~\ref{e7zh} (right).

\medskip
The singularity $E_8$ has two perturbations, whose zero sets look topologically as in Fig.~\ref{e8pp}. The first of them is explicitly constructed (up to the sign) in page 17 of \cite{AC}, and the other one is realized by the same algorithm if before the last contraction we shift a triple intersection point to the different side.
\unitlength 0.7mm
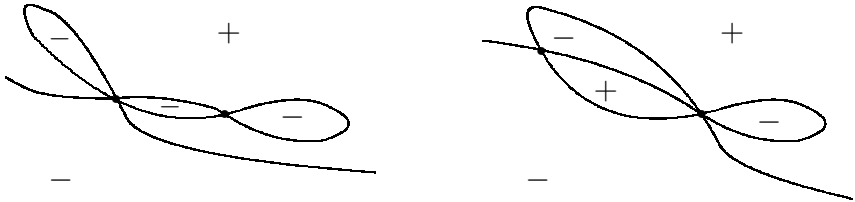
\begin{figure}
{
\begin{picture}(65,36)
\bezier{70}(0,23)(5,20)(7,20)
\bezier{70}(7,20)(11,19)(20.2,19)
\bezier{100}(20.2,19)(30,20)(40,17)
\bezier{100}(40,17)(50,10)(60,11)
\bezier{70}(60,11)(70,14)(60,18)
\bezier{100}(60,18)(55,20)(45,17)
\bezier{150}(45,17)(25,10)(5,31)
\bezier{100}(5,31)(1,39)(7,36)
\bezier{200}(7,36)(13,35)(23,15)
\bezier{200}(23,15)(28,8)(70,5)
\put(20.9,18.9){\circle*{1.5}}
\put(41.5,16){\circle*{1.5}}
\put(8,2){$-$}
\put(8,29){\small $-$}
\put(29,16){\small $-$}
\put(52,14){$-$}
\put(40,30){$+$}
\end{picture} \qquad \qquad
\begin{picture}(70,36)
\bezier{200}(0,30)(25,27)(40,17)
\bezier{100}(40,17)(50,10)(60,11)
\bezier{70}(60,11)(70,14)(60,18)
\bezier{100}(60,18)(55,20)(45,17)
\bezier{150}(45,17)(20,10)(10,30)
\bezier{100}(10,30)(6,38)(12,36)
\bezier{200}(12,36)(35,30)(45,10)
\bezier{150}(45,10)(48,5)(70,0)
\put(11.1,28.1){\circle*{1.5}}
\put(41.5,16){\circle*{1.5}}
\put(8,2){$-$}
\put(13,29){$-$}
\put(21,19){$+$}
\put(52,13){$-$}
\put(45,30){$+$}
\end{picture}
}
\caption{Degenerate perturbations for $E_8$}
\label{e8pp}
\end{figure}

The construction of \S \ref{basic} applied to these two perturbations gives us 1-cycles in components with passports $(3, 3, 0)$ and $(2, 3, 1)$. The involution $f_\lambda(x, y) \mapsto -f_\lambda(-x, -y)$ of the space of deformation (\ref{e8vd}) turns these cycles to other two, which lie in components with passports $(0, 3, 3)$ and $(1, 3, 2)$. 

The perturbation $x^3 + y^5 + \varepsilon y^4$ splits the singularity $E_8$ into two real critical points of types $+E_6$ and $A_2$. We can further perturb these two points independently, removing the point of type $A_2$ from the real domain and moving the point $+E_6$ as shown in Fig.~\ref{e6proof} (left); so we get a cycle in the component with passport $(2, 2, 0)$. The symmetric perturbation gives us passport $(0, 2, 2)$. 

The perturbation $x^3 + y^5 + \varepsilon x^2 y$ has a singularity of type $D_6^+$ and no other real critical points. We have seen in page \pageref{1d2p} that in a neighborhood of such a function there is a non-trivial 1-cycle in the component with passport $(1, 2, 1)$. So, we have realized 1-cycles in components with seven different passports; the winding number takes value 1 on all these cycles.

\section{Two-dimensional homology of complements of caustics of some non-simple singularities}
\label{H2}

\subsection{$J_{10}^3$} The singularities of class $J_{10}^3$ have the normal form 
\begin{equation}
x(x+y^2)\left(x-\alpha y^2\right), \quad \alpha >0, 
\label{j103}
\end{equation} 
and the Milnor number equal to 10. A miniversal deformation of any function (\ref{j103}) can be chosen in the form
\begin{equation}
f_0 + \lambda_1 + \lambda_2 x + \lambda_3 y + \lambda_4 x y + \lambda_5 y^2 + \lambda_6 x y^2 + \lambda_7 y^3 + \lambda_8 x y^3 + \lambda_9 y^4 + \lambda_{10} x y^4. 
\label{j103vd}
\end{equation}
For arbitrary $\varepsilon >0$, this deformation contains the perturbation
\begin{equation}
\tilde f = x (x + y^2 - \varepsilon^2)(x- \alpha (y^2-\varepsilon^2))
\end{equation} 
 of the function (\ref{j103}). This perturbation has a Morse minimum, a Morse maximum, and two critical points of type $D_4^-$; its zero level set looks as \ \ 
\unitlength 1mm
\begin{picture}(22,7)
\put(0,2){\line(1,0){20}}
\bezier{150}(0,-1)(10,9)(20,-1)
\bezier{150}(0,7)(10,-9)(20,7)
\put(3.7,2){\circle*{1}}
\put(16.3,2){\circle*{1}}
\end{picture} (where the line $\{x=0\}$ is horizontal). 

In a neighborhood of the point $\tilde f$, the deformation (\ref{j103vd}) provides a versal deformation of corresponding multisingularity, therefore the pair $( {\mathbb R}^{10}, \mbox{the caustic})$ is locally diffeomorphic in such a neighborhood to the product of $ {\mathbb R}^4$ and the pair $( {\mathbb R}^3 \times  {\mathbb R}^3, C),$ where each of two copies of $ {\mathbb R}^3$ is identified with the parameter space of the shortened versal deformation of type $D_4^-$ (see Fig.~\ref{pyramid}), and the subset $C \subset  {\mathbb R}^3 \times  {\mathbb R}^3$ consists of all points whose projection to at least one of factors $ {\mathbb R}^3$ belongs to the caustic of this deformation. The product of 1-cycles embracing the pyramids in these spaces $ {\mathbb R}^3$ (see \S \ref{basic}) is a torus $T^2$ in the space $ {\mathbb R}^3 \times  {\mathbb R}^3 \setminus C$, hence it defines also a torus in the complement of the caustic of our singularity. This torus lies in a component of this complement with passport $(1, 4, 1)$.

\begin{theorem}
The map of entire homology group $H_*(T^2, {\mathbb Z})$ of the torus described in the previous paragraph to the homology group of the complement of the caustic of singularity $J_{10}^3$ induced by identical embedding is monomorphic.
\end{theorem}

\noindent
{\it Proof.} The restriction of the derivative $\frac{\partial f_\lambda}{\partial x}$ of any perturbation $f_\lambda$ of the form (\ref{j103vd}) to any line $\{y = \mbox{const}\}$ is a polynomial of degree 2 with non-zero leading term, so $f_\lambda$ cannot have more than two critical points in such a line.
 Therefore, the map of our component of the complement of the caustic to $B( {\mathbb R}^2, 4)$ defined in \S \ref{bascoc} acts in fact into the subset $\hat B( {\mathbb R}^2,4) \subset B( {\mathbb R}^2, 4)$ consisting of configurations, none three points of which lie on the same vertical line. The cohomology group of this subset follows immediately from the calculations of \cite{vain}, as $\hat B( {\mathbb R}^2,4)$ is a union of several cells of the cell structure considered in \cite{fuks}, \cite{vain}. Namely, $\hat B( {\mathbb R}^2,4)$ consists of unique 8-dimensional cell
 \unitlength 0.5mm
\begin{picture}(21,10)
\put(10.5,3){\oval(21,10)}
\put(3,3){\circle*{2}} 
\put(8,3){\circle*{2}} 
\put(13,3){\circle*{2}} 
\put(18,3){\circle*{2}} 
\end{picture}
consisting of configurations, all four points of which have different $x$-coordinate, three 7-dimensional cells 
\begin{picture}(16,10)
\put(8,3){\oval(16,10)}
\put(3,1){\circle*{2}} 
\put(3,5){\circle*{2}} 
\put(8,3){\circle*{2}} 
\put(13,3){\circle*{2}} 
\end{picture} , 
\begin{picture}(16,10)
\put(8,3){\oval(16,10)}
\put(3,3){\circle*{2}} 
\put(8,1){\circle*{2}} 
\put(8,5){\circle*{2}} 
\put(13,3){\circle*{2}} 
\end{picture} and 
\begin{picture}(16,10)
\put(8,3){\oval(16,10)}
\put(3,3){\circle*{2}} 
\put(8,3){\circle*{2}} 
\put(13,1){\circle*{2}} 
\put(13,5){\circle*{2}} 
\end{picture} of 4-configurations with only three different values of $x$-coordinate, and one 6-dimensional cell 
\begin{picture}(11,10)
\put(5.5,3){\oval(11,10)}
\put(3,1){\circle*{2}} 
\put(3,4){\circle*{2}} 
\put(8,2){\circle*{2}} 
\put(8,5){\circle*{2}} 
\end{picture}. 
According to \cite{vain}, all incidence coefficients of these cells are trivial, so by Poincar\'e duality we have $H^1(\hat B( {\mathbb R}^2,4),  {\mathbb Z}) \simeq  {\mathbb Z}^3$ and $H^2(\hat B( {\mathbb R}^2,4),  {\mathbb Z}) \simeq  {\mathbb Z}.$ It follows immediately from the construction that the intersection index of the image of our torus with the cell \begin{picture}(11,10)
\put(5.5,3){\oval(11,10)}
\put(3,1){\circle*{2}} 
\put(3,4){\circle*{2}} 
\put(8,2){\circle*{2}} 
\put(8,5){\circle*{2}} 
\end{picture} is equal to $1$ or $-1$ depending on the choice of orientations, and intersection indices of images of two 1-dimensional generators of the torus with 7-dimensional cells are equal to $(\pm 1, 0, 0)$ and $(0, 0, \pm 1)$. \hfill $\Box$

\subsection{Class $P_8^2$} 

The singularity classes $P_8^1$ and $P_8^2$ consist of functions of corank 3 with non-degenerate cubic part. Any singularity of this type depending on three variables can be reduced by a local diffeomorphism of the argument space to a non-degenerate homogeneous polynomial of degree 3. Its zero set defines in $ {\mathbb R}P^2$ a cubic curve with one or two components, in the correspondence with the upper index of the name of the class. The Milnor numbers of all singularities of this these classes are equal to 8, and their miniversal deformations can be chosen in the form
\begin{equation}
\label{vdp82}
f + \lambda_1 + \lambda_2 x + \lambda_3 y + \lambda_4 z + \lambda_5 x y + \lambda_6 x z + \lambda_7 y z + \lambda_8 x y z \ ,
\end{equation}
where $f$ is the initial homogeneous polynomial.

\begin{theorem}
\label{th44}
For any sufficiently small $\varkappa > 0,$
the 2-dimensional homology group of the complement of the caustic of versal deformation (\ref{vdp82}) of function 
\begin{equation}
f_\varkappa = (x-z)(x+z)(\varkappa x-z) - y^2 z
\label{p82}
\end{equation}
 of class $P_8^2$ is non-trivial. 
\end{theorem}

\unitlength 1.3 mm
 \begin{figure}
\begin{picture}(80,43)
\bezier{100}(20,10)(30,10)(30,20)
\bezier{100}(30,20)(30,30)(20,30)
\bezier{100}(20,30)(10,30)(10,20)
\bezier{100}(10,20)(10,10)(20,10)
\bezier{150}(18.,3)(35.,3)(35.,20)
\bezier{150}(35.,20)(35.,37)(18.,37)
\bezier{150}(18.,37)(1.,37)(1.,20)
\bezier{150}(1.,20)(1.,3)(18.,3)
\put(18,20){\circle{0.9}}
\put(0,20){\vector(1,0){80}}
\put(20,0){\vector(0,1){43}}
\put(10,20){\circle*{0.7}}
\put(30,20){\circle*{0.7}}
\put(65,20){\circle*{0.7}}
\put(20,20){\circle*{0.7}}
\put(18,20){\vector(-3,4){10.7}}
\put(47,30){$-$}
\put(75,30){$+$}
\put(22,24){$+$}
\put(22,41){$y$}
\put(77,17){$x$}
\put(62,15.7){$\frac{1}{\varkappa}$}
\put(4.7,16.5){\small $-1$}
\put(30.5,16.5){\small $1$}
\put(20.3,16.5){\small $0$}
\bezier{80}(65,20)(65,26)(68,30)
\bezier{80}(65,20)(65,14)(68,10)
\bezier{80}(68,30)(71,36)(72,40)
\bezier{80}(68,10)(71,4)(72,0)
\end{picture}
\caption{Affine chart \ $z=1$}
\label{pc}
\end{figure}
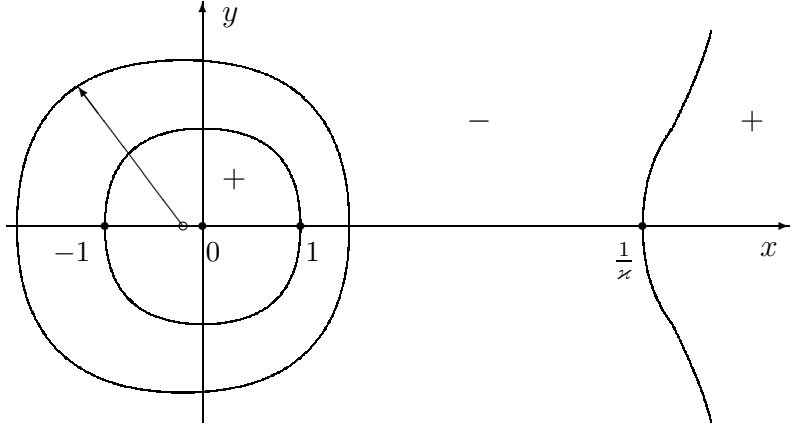

To formulate a more concrete statement, consider the following $n+1$ local systems $\pm  {\mathbb Z}_{\{i\}}$, $i=0, 1, \dots, n$,
on any connected component of the complement of the caustic of a function in $n$ variables. Each of these systems is locally isomorphic to $ {\mathbb Z}$, but the loops in our component, defining odd permutations of critical points with Morse index $i$, act on the fibers as multiplication by $-1$.
\medskip 

\noindent 
{\bf Theorem $\mbox{\ref{th44}}'$}. The complement of the caustic of singularity (\ref{p82}) with sufficiently small $\varkappa>0$ 
has a connected component with passport $(0, 2, 2, 0)$ such that the 2-dimensional homology group of this component with coefficients in $ {\mathbb Z}_2$ is non-trivial, and the 2-dimensional homology group with coefficients in the local system $\pm  {\mathbb Z}_{\{1\}}$ contains a subgroup isomorphic to $ {\mathbb Z}$.
\medskip

\noindent
The proof of this theorem takes the rest of this section. The non-trivial homology classes proving Theorem $\mbox{\ref{th44}}'$ will be realized by the fundamental classes of the Klein bottle $K^2$ embedded somehow into the complement of the caustic. 
\medskip

\noindent
{\bf Remark.}
We will consider also the degenerate function $f_0$ and its deformation given respectively by formulas (\ref{p82}) and (\ref{vdp82}) with $\varkappa=0$. Although in this case this deformation is not versal, all our constructions in the real domain will hold for it as well. 
\medskip

The polynomial (\ref{p82}) with arbitrary $\varkappa \in [0,1)$ is {\em hyperbolic} with respect to the line $\{x=y=0\}$, see e.g. \cite{APLT}, \S IV.2. Its zero set consists of two components, one of which is homeomorphic to a quadratic cone, and the other to a hyperplane in $ {\mathbb R}^3$. The corresponding sets in the affine chart $\{z=1\}$ of the projectivization $ {\mathbb R}P^2$ of $ {\mathbb R}^3$ are shown respectively by the interior oval on the left of Fig.~\ref{pc} and by the non-closed curve on the right. (If $\varkappa=0$ then this curve becomes a straight line.)
We will assume that the argument space $ {\mathbb R}^3$ is Euclidean with the standard metric defined by coordinates $x, y, z$. Recall that the {\it norm} of a quadratic function $\|L\|$ is the maximal value of $|L|$ on the unit ball in $ {\mathbb R}^3$.

\begin{lemma}
\label{lem98}
1. If $\varkappa \geq 0$ is small enough, and $L:  {\mathbb R}^3 \to  {\mathbb R}$ is an arbitrary 
quadratic function of rank 1 taking zero value on the line $\{x=y=0\}$, then the function $f_\varkappa + L$ has a singularity of type $D_4^-$ at the origin, and additionally exactly two Morse real critical points. Morse indices of the last two points are equal to 2 $($respectively, 1$)$ if $L$ is non-negative $($respectively, non-positive$)$. 

2. There are positive constants $\theta< \Theta$ such that the distances of both these critical points from the origin in $ {\mathbb R}^3$ and from the plane $\ker L$ separating them belong to the segment $\left[\theta \|L\|, \Theta \|L\|\right]$ for any such function $L$. 
\end{lemma}

\noindent
{\it Proof of Lemma \ref{lem98}.} 
 It is enough to prove statement 1 of this lemma only for functions $f_\varkappa + L$ such that the quadratic function $L$ is non-negative with norm 1, i.e. has the form $( \alpha x + \beta y)^2$ where $\alpha^2 + \beta^2=1$. Indeed, for any such function $L$ and any non-zero constant $t$ we have 
\begin{equation}
\label{scaling}
(f_\varkappa + t L)(t x, t y, t z) \equiv t^3 (f_\varkappa + L)(x, y, z).
\end{equation}

Denote by $\bigcirc$ the set of quadratic functions $L$ of this type (i.e. of rank 1, norm 1, and vanishing on the $z$ axis). This set is obviously diffeomorphic to a circle.

Statement 2 of Lemma \ref{lem98} follows from statement 2. Indeed, for any $L \in \bigcirc$ by (\ref{scaling}) the discussed distances of critical points of functions $f_\varkappa + t L$ are proportional to $t$, so it is enough to prove this statement for $L$ from the compact set $\bigcirc$. The distances in question define positive functions on this set, and hence are separated by constants from 0 and infinity.

The corank of function $f_\varkappa + L$ is equal to 2, and the restriction of this function to the two-dimensional kernel \ $\ker L$ \ of its quadratic part is a homogeneous polynomial of degree 3 vanishing on three real lines. Therefore in some linear coordinates $(\xi, \eta, \zeta)$ in $ {\mathbb R}^3$ the polynomial $f_\varkappa$ consists of the principal part $\xi^2 + \eta^2 \zeta - \zeta^3$ and several monomials of higher quasihomogeneous degree with respect to weights $(3 : 2 : 2)$ of these coordinates. The standard algorithm of the reduction to normal forms (see e.g. \cite{AVGZ}) removes all these additional monomials of its Taylor expansion at the origin and reduces the germ of our function at the origin to the standard normal form of class $D_4^-$.

At all the other critical points of the function $f_\varkappa + L$, the partial derivative $\frac{\partial f_\varkappa}{\partial z}$ is equal to 0 (as the gradient of $f_\varkappa$ should be opposite to that of $L$).
The set of points in $ {\mathbb R}^3$ satisfying this condition is the cone 
\begin{equation}
(x + \varkappa z)^2+y^2 - (3 + \varkappa^2)z^2=0.
\label{cone}
\end{equation} In the affine chart in $ {\mathbb R}P^2 $ defined by the plane with equation $z=1$ in $ {\mathbb R}^3$ it is represented by the circle of radius $\sqrt{3+ \varkappa^2}$ with the center at the point $(-\varkappa, 0)$; see the greater oval in Fig.~\ref{pc}. In particular, this cone does not intersect the set $\{f_\varkappa =0\}$ in $ {\mathbb R}^3 \setminus 0$ and separates all its components. For any non-zero point $A$ of the cone (\ref{cone}), there is only one quadratic function $L(A)$ of rank 1 vanishing on the line $\{x=y=0\}$ such that $f_\varkappa+L$ has a critical point at $A$. If $A$ belongs to the half of this cone where $z>0$ (respectively, $z<0$), then this function $L(A)$ is everywhere non-negative (respectively, non-positive). Indeed, in this half a small shift along the vector $\mbox{grad} f_\varkappa$ decreases (respectively, increases) the distance from the line $\{x=y=0\}$, and the vector $\mbox{grad} \, L(A)$ should be opposite to it.

The last part of the first assertion of Lemma \ref{lem98} means that this map $\{ A \mapsto L(A)\}$ is a two-fold covering of the smooth part of cone (\ref{cone}) over the set of all quadratic functions of rank 1 vanishing on the line $\{x=y=0\}$. Formula (\ref{scaling}) allows us to reformulate this statement in terms of compact quotients of these sets by the scalings. Namely, denote by $S^1_\varkappa$ the greater oval of Fig.~\ref{pc}, 
i.e. the intersection of the cone (\ref{cone}) and the affine plane $\{z=1\}$. 
Denote by $\Phi_\varkappa$ the map $S_\varkappa^1 \to \bigcirc$, which sends any point $A \in S^1_\varkappa$ to the unique function $L \in \bigcirc$ such that one of the functions $f_\varkappa + t L$, $t >0$, has a critical point at $A$. It remains to prove that for sufficiently small $\varkappa \geq 0$ this map $\Phi_\varkappa$ is a smooth two-fold covering, and all these critical points are Morse with Morse indices equal to 2. 

To do it, notice that the analogous statement (and entire Lemma \ref{lem98}) is true for the degenerate function (\ref{p82}) with $\varkappa=0.$ Indeed, both this function $f_0$ and the circle $S^1_0$ are invariant under rotations around the axis $x=y=1$, and the map $\Phi_0$ commutes with these rotations. Therefore it is sufficient to check this statement for the model linear function $L \equiv x^2$, which is elementary. 

For all sufficiently small $\varkappa \geq 0$ define the maps $\rho_\varkappa : S^1_0 \to S^1_\varkappa$ as follows. For any point $A \in S^1_0$ consider the point $\Phi_0(A) \in \bigcirc$, i.e. 
the unique quadratic function vanishing on the line $\{x=y=0\}$ such that $A$ is a critical point of the function $f_0 + \Phi_0(A)$. Then,
as we will show in the next paragraph, the function $f_\varkappa + \Phi_0(A)$ has only one critical point neighboring to $A$ in $ {\mathbb R}^3$. We denote this point by 
$\tilde A_\varkappa (A) \in  {\mathbb R}^3$ and define 
 $\rho_\varkappa(A) \in S^1_\varkappa$ as the point of the plane $\{z=1\}$ proportional to $\tilde A_\varkappa (A)$ in the vector space $ {\mathbb R}^3$. 

This map $\rho_\varkappa$ is well-defined and is a diffeomorphism for sufficiently small $\varkappa \geq 0$; for $\varkappa=0$ it is the identical map. Indeed, consider the product of spaces $ {\mathbb R}^1 \times S^1_0 \times  {\mathbb R}^3$ whose elements are called $\varkappa, A$ and $(x, y, z)$, and three conditions \begin{equation}
\label{pard}
\frac{\partial}{\partial x} (f_\varkappa + \Phi_0(A)) = \frac{\partial}{\partial y} (f_\varkappa + \Phi_0(A)) = \frac{\partial}{\partial z} (f_\varkappa + \Phi_0(A)) = 0.
\end{equation}
The matrix of partial derivatives of three functions (\ref{pard}) upon $x, y, $ and $z$ is the Hessian of the function $f_\varkappa + \Phi_0(A)$. It is easy to check that for $\varkappa =0$ and arbitrary $A \in S^1_0$ its determinant is equal to 3, therefore by implicit function theorem the map $(\varkappa, A) \mapsto \tilde A_\varkappa(A)$ is well-defined and smooth in a neighborhood of the circle $0 \times S^1_0$ in $ {\mathbb R}^1 \times S^1_0$. The restriction of this map to this circle is a smooth (identical) embedding, hence its restrictions to all neighboring circles $\varkappa \times S^1_0$ also are smooth embeddings. Moreover, the radial projection of the image of this circle from $ {\mathbb R}^3$ 
to the plane $\{z=1\}$ has no singular points, hence the same is true for the neighboring embedded circles, and our maps $\{ A \mapsto \tilde A_\varkappa(A)\}$ for all small $\varkappa$ are smooth immersions of $S^1_0 $ to $ {\mathbb R}^2$, the images of which coincide with corresponding circles $S^1_\varkappa$.

We get the commutative diagram
\begin{eqnarray*}
S^1_0 \ \ \ \ & \stackrel{\rho_\varkappa}{\longrightarrow} & \ \ \ \ S^1_\varkappa \\
 \Phi_0 \searrow & & \swarrow \Phi_\varkappa \\
 & \bigcirc & 
\end{eqnarray*}
in which $\Phi_0$ is a two-fold covering and $\rho_\varkappa$ is a diffeomorphism, hence $\Phi_\varkappa$ also is a covering. \hfill $\Box$ 
\medskip

\noindent
{\bf Remark.}
The collection of all {\em complex} critical points of functions $ f_\varkappa + L$ with small $\varkappa \neq 0$ {\em is not} a small perturbation of that for the model function $f_0 + L$ : it additionally contains two critical points in the complex domain.
\medskip

Let again $L$ be a quadratic function of class \ $\bigcirc$, \ and let $(\xi, \eta, \zeta)$ be linear coordinates in $ {\mathbb R}^3$ such that $L \equiv \xi^2$, $\eta $ and $\zeta$ are constant on lines orthogonal to \ $\ker L$, \ and the restriction of $f_\varkappa$ to \ $\ker L$ \ is equal to $3\eta \zeta - \zeta^3$, cf. \S \ref{basic}. For an arbitrary small number $\varepsilon >0$ consider the family of perturbations of the function $f_\varkappa$
depending on the parameter $\tau \in [0, 2\pi]$ and defined by formula 
\begin{equation}
f_\varkappa + \varepsilon L - 3 \varepsilon^4 (\sin \tau \cdot \eta + \cos \tau \cdot \zeta),
\label{kb}
\end{equation}
cf. (\ref{param}). Functions of this family form a subset diffeomorphic to $S^1$ in the space of deformation (\ref{vdp82}). This subset does not depend on the choice of linear coordinates $(\xi, \eta, \zeta)$ as above, because all such choices differ from each other by orthogonal transformations moving this family to itself. Such subsets defined by formula (\ref{kb}) for all choices of $L \in \bigcirc$ \ sweep out a surface with the structure of a fiber bundle over \ $\bigcirc$ \ with fiber $S^1$. Moving $L$ once along the loop \ $\bigcirc$ \ we rotate the corresponding plane \ $\ker L$ \ by angle $\pi$ around a fixed axis and hence change the orientation of this plane and of the fiber of our fiber bundle; so the space of this bundle is homeomorphic to the Klein bottle. 

\begin{lemma} If $\varepsilon >0$ is sufficiently small, then

1. Any function of the form $($\ref{kb}$)$ has exactly four real critical points, all of which are Morse; in particular the family $($\ref{kb}$)$ lies in the complement of the caustic.

2. Two of these points lie on distance of order \ $\varepsilon$ \ from the origin and from the plane \ $\ker L$, and are separated by this plane; their Morse indices are equal to 2. 

3. Two other critical points lie on distance of order \ $\varepsilon^2$ \ from the origin and from one another, and their distances from the plane \ $\ker L$ \ are estimated from above by $ \varepsilon^{8/3}$; the Morse indices of last two critical points are equal to $1$.
\label{lem3}
\end{lemma}

\noindent
{\it Proof.} By Lemma \ref{lem98}, the function $f_\varkappa + \varepsilon L$ has two Morse critical points on distance of order $\varepsilon$ from the origin and from the plane $\ker L$, separating these two points. 
By scaling argumentation (cf. (\ref{scaling})) the norms of inverse Hessians of this function in a neighborhood of the union of these points over all $L \in \bigcirc$ are estimated from below by the number $c \varepsilon,$ where $c$ is a positive constant not depending on $L$; therefore the further minor perturbation by subtracting
$3\varepsilon^4 (\sin \tau \cdot \eta + \cos \tau \cdot \zeta)$ does not move these points seriously and does not change their Morse indices (which by Lemma \ref{lem98} are equal to 2). 

On the other hand, this perturbation splits the singularity of type $D_4^-$ of the function $f_\varkappa + \varepsilon L$ at the origin into two real points with Morse index 1 and two non-real points. Indeed, in coordinates $(\xi, \eta, \zeta)$ the function $f_\varkappa + \varepsilon L$ equals
\begin{equation}
3 \eta^2 \zeta - \zeta^3 + \varepsilon \xi^2 + \xi \omega_2(\eta, \zeta) + \xi^2 \omega_1(\eta, \zeta) + \xi^3 \omega_0,
\label{nf}
\end{equation}
where $\omega_k$ is a homogeneous polynomial of degree $k$. 
Let us break the function $f_\varkappa + \varepsilon L - 3 \varepsilon^4 (\sin \tau \cdot \eta + \cos \tau \cdot \zeta)$ 
into two parts, $$P_\varkappa \equiv 3\eta^2 \zeta - \zeta^3 + \varepsilon \xi^2- 3 \varepsilon^4 (\sin \tau \cdot \eta + \cos \tau \cdot \zeta) \quad \mbox{and} \quad Q_\varkappa \equiv \xi \omega_2(\eta, \zeta) + \xi^2 \omega_1(\eta, \zeta) + \xi^3 \omega_0.$$

By \S \ref{basic} the function $P_\varkappa $ has exactly four Morse critical points in $ {\mathbb C}^3$ with coordinates 
\begin{equation}
(\xi, \eta, \zeta) = \pm~\varepsilon^2\left(0, \, \cos \frac{\tau}{2}, \, \sin \frac{\tau}{2}\right) \qquad \mbox{and} \qquad \pm~i \varepsilon^2 \left( 0, \, -\sin \frac{\tau}{2}, \, \cos \frac{\tau}{2} \right);
\label{appr}
\end{equation} 
two first of them are real with Morse index 1, ant two others are imaginary.

Let us define small neighborhoods $U_j \subset  {\mathbb C}^3$, $j=1, 2, 3, 4$, of all these four points $(0, \eta_j, \zeta_j)$ by respective conditions 
\begin{equation}
|\eta-\eta_j|^2+|\zeta-\zeta_j|^2 \leq \varepsilon^{14/3} 
\label{etbo}
\end{equation}
and by the common inequality \begin{equation}
 |\xi| \leq \varepsilon^{8/3} .
\label{xibo}
\end{equation} 

\begin{proposition} 
\label{pro1z}
If $\varepsilon$ is small enough, then for each $j=1, \dots, 4$
\begin{equation}
\label{estsid}
\left|\frac{\partial P_\varkappa }{\partial \eta}\right|^2 + \left|\frac{\partial P_\varkappa }{\partial \zeta}\right|^2 > \left|\frac{\partial Q_\varkappa}{\partial \eta}\right|^2 + \left|\frac{\partial Q_\varkappa}{\partial \zeta}\right|^2 
\end{equation}
on the part of the boundary of the domain $U_j$, where the corresponding inequality $($\ref{xibo}$)$ becomes an equality, and 
\begin{equation}
\left|\frac{\partial P_\varkappa }{\partial \xi}\right| > \left|\frac{\partial Q_\varkappa}{\partial \xi}\right|
\label{estup}
\end{equation} 
on the part of its boundary where the inequality $($\ref{etbo}$)$ becomes an equality. 
\end{proposition} 

\noindent
{\it Proof of Proposition \ref{pro1z}}. If $|\xi|=\varepsilon^{8/3}$ then $\left|\frac{\partial P_\varkappa}{\partial \xi}\right| = 2\varepsilon^{11/3}$, while three summands of $\left|\frac{\partial Q_\varkappa}{\partial \xi}\right|$ have asymptotic behaviors of types $O(\varepsilon^4), O(\varepsilon^{14/3}) $ and $O(\varepsilon^{16/3})$ respectively when $\varepsilon$ tends to zero. 

The shift of the origin of our coordinate system \ $(\xi, \eta, \zeta)$ \ to either of four critical points (\ref{appr}) of $P_\varkappa $ transforms the function $P_\varkappa - \varepsilon \xi^2$ (i.e. the part of $P_\varkappa$ depending on variables $\eta$ and $\zeta$) to a polynomial of degree 3 with zero linear part. The norm of the gradient of the quadratic part of this polynomial at the piece of the boundary of this critical point, where (\ref{etbo}) becomes an equality, is equal identically to $6 \varepsilon^{13/3}$.

The gradient of the cubic part of this function has asymptotic type $O(\varepsilon^{14/3})$, and the asymptotic types of partial derivatives upon variables $\eta$ and $\zeta$ of the summand $\varepsilon \xi^2$ of $P_\varkappa$ and of three summands of $Q_\varkappa$ are respectively $0$, $O(\varepsilon^{8/3 +2}) = O(\varepsilon^{14/3}), O(\varepsilon^{16/3})$ and $0$. \hfill $\Box$ \medskip

\noindent
{\it End of the proof of Lemma \ref{lem3}.}
By Proposition \ref{pro1z} each function from the one-parametric family $P_\varkappa + t Q_\varkappa, $ $t \in [0,1]$, connecting the functions $P_\varkappa$ and 
(\ref{kb}), has exactly one complex critical point in any of four neighborhoods $U_j$. 
All these critical points are Morse; those of them which are real for $t=0$ remain real and preserve their Morse indices for all values of $t$. 
\hfill $\Box$ \medskip

So, the component of the complement of the caustic of the function (\ref{p82}) with sufficiently small $\varkappa \geq 0$ that contains all functions $($\ref{kb}$)$ consists of morsifications with exactly two real critical points with Morse index 1 and two critical points with index 2. 
Denote by $G_\varkappa$ the map from this component to $ {\mathbb R}P^2$ sending its point $\lambda$ to the direction of the line passing through two critical points of the corresponding function $f_\lambda$ with Morse index 1. By definition, a loop in this component goes into the non-trivial element of $\pi_1( {\mathbb R}P^2)$ if and only if it permutes these two critical points or, equivalently, it changes the orientation of the local system $\pm  {\mathbb Z}_{\{1\}}$. 

\begin{proposition}
\label{pro3z}
The restriction of the map $G_\varkappa$ to the embedded Klein bottle consisting of functions $($\ref{kb}$)$ maps the fundamental class of $K^2$ with coefficients in the orientation sheaf to that of $ {\mathbb R}P^2$.
\end{proposition}

\noindent 
{\it Proof.} 
Denote by $\tilde K^2$ the set of all pairs (a 2-dimensional subspace in $ {\mathbb R}^3$ passing through the axis $\{x=y=0\}$, a 1-dimensional subspace of it). This set is obviously diffeomorphic to the Klein bottle. Forgetting the first elements of these pairs maps it to $ {\mathbb R}P^2$; this map is a diffeomorphism over all points of $ {\mathbb R}P^2$ but only one, and defines an isomorphism between 2-homology groups of $\tilde K^2$ and $ {\mathbb R}P^2$ with coefficients in orientation sheaves. 

Further, let $K^2(0)$ be the set of all functions of the form (\ref{kb}) with $\varkappa=0$ in the space of Morse perturbations of the degenerate function $f_0$. Let $\rho: K^2(0) \to \tilde K^2$ be the map associating with any function of this form the pair consisting of the corresponding plane $\ker L$ and the line passing through the real critical points of the restriction of this function to this plane $\ker L$. This map is a diffeomorphism by the construction and \S~\ref{basic}. The composition $K^2(0) \to \tilde K^2 \to  {\mathbb R}P^2$ of this map and the tautological map considered in the previous paragraph is homotopic to the restriction of the map $G_0$ to $K^2(0)$. Indeed, by Lemma \ref{lem3} and Proposition \ref{pro1z} the images of any two points under these two maps lie on distance of order $\varepsilon^{2/3}$ from one another in the standard metric on $ {\mathbb R}P^2$ induced from the unit sphere in $ {\mathbb R}^3$. Therefore the statement of Proposition \ref{pro3z} is true for map $G_0$. Finally, the entire our construction of embedded Klein bottles and maps $G_\varkappa$ for sufficiently small positive $\varkappa$ is a small perturbation of that for $\varkappa =0$. 
This finishes the proof of Theorem $\mbox{\ref{th44}}'$ for the case of $\pm  {\mathbb Z}_{\{1\}}$-coefficients. The statement of this theorem concerning constant $ {\mathbb Z}_2$-coefficients is its reduction mod 2. \hfill $\Box$ 

\subsection{$P_8^1$ conjecture} The class $P_8^1$ is defined similarly to $P_8^2$, but with only one component of zero set in $ {\mathbb R}^3$, homeomorphic to $ {\mathbb R}^2$. Its easiest representative is the function $x^3 + y^3 + z^3$. Consider the perturbation $x^3 - 3 \varepsilon^2 x + y^3 - 3 \varepsilon^2 y + z^3$ of this function. It has four real critical points of type $A_2$ with coordinates $(\pm \varepsilon, \pm \varepsilon, 0)$ (where the signs $\pm$ are independent). We can perturb these points independently within the versal deformation. Let us do it so that the points with coordinates $(\varepsilon, \varepsilon)$ and $(-\varepsilon, -\varepsilon)$ split into imaginary critical points, and each of the remaining two points splits into a pair of real Morse points. The passport of the resulting function is $(0, 2, 2, 0)$. 

\begin{conjecture}
The two-dimensional cohomology group with coefficients in $ {\mathbb Z}_2$ of the component of the complement of the caustic of the morsification just constructed is non-trivial. Namely, the 1-cohomology classes watching the permutations of two critical points of the same Morse index are non-trivial, as well as their cohomological product.
\end{conjecture}

\end{document}